\documentclass[titlepage,11pt]{article}

\usepackage{amsmath, amssymb}
\numberwithin{equation}{section}

\begin{document}
\baselineskip 1.65em

\newcommand{\re}{\mathop{{\rm Re}}}
\newcommand{\im}{\mathop{{\rm Im}}}
\newcommand{\interior}{\mathop{{\rm int}}}
\newcommand{\rad}{\mathop{{\rm rad}}}
\newcommand{\alp}{\alpha}
\newcommand{\bet}{\beta}
\newcommand{\gam}{\gamma}
\newcommand{\Gam}{\Gamma}
\newcommand{\del}{\delta}
\newcommand{\eps}{\varepsilon}
\newcommand{\zet}{\zeta}
\newcommand{\lam}{\lambda}
\newcommand{\Lam}{\Lambda}
\newcommand{\vph}{\varphi}
\newcommand{\sig}{\sigma}
\newcommand{\ome}{\omega}
\newcommand{\Ome}{\Omega}
\newcommand{\IR}{{\bf R}}
\newcommand{\IC}{{\bf C}}
\newcommand{\ZZ}{{\bf Z}}
\newcommand{\wt}{\widetilde}
\newcommand{\wh}{\widehat}
\newcommand{\ds}{\displaystyle}
\newcommand{\ts}{\textstyle}
\newcommand{\bxx}{\begin{it}}
\newcommand{\exx}{\end{it}}
\newcommand{\xx}{\it}
\newcommand{\vv}{\vspace{.2in}}
\newcommand{\bb}{\hspace*{-.08in}}
\newcommand{\fp}{\hspace{1mm}\rule{2mm}{2mm}}
\newcommand{\be}[1]{\begin{equation}}
\newcommand{\ee}[1]{\label{#1}\end{equation}}

\title{Analyticity~and~Nonanalyticity~of Solutions~of~Delay-Differential~Equations}

\author{John Mallet-Paret\thanks{Partially supported by The Center for Nonlinear Analysis at Rutgers University.}
\\Division of Applied Mathematics
\\Brown University
\\Providence, Rhode Island 02912
\\{\tt john\_mallet-paret@brown.edu}
\\{ }
\and
Roger D. Nussbaum\thanks{Partially supported by NSF Grant DMS-1201328 and by The Lefschetz Center for Dynamical Systems
at Brown University.}
\\Department of Mathematics
\\Rutgers University
\\Piscataway, NJ 08854
\\{\tt nussbaum@math.rutgers.edu}
\\{ }}
\date{April 30, 2013}

\maketitle

\begin{abstract}
We consider the equation
$$
\dot x(t)=f(t,x(t),x(\eta(t)))
$$
with a variable time-shift $\eta(t)$.
Both the nonlinearity $f$ and the shift function $\eta$ are
given, and are assumed to be analytic (that is, holomorphic)
functions of their arguments. Typically
the time-shift represents a delay, namely that $\eta(t)=t-r(t)$
with $r(t)\ge 0$.
The main problem considered is to determine when
solutions (generally $C^\infty$ and often periodic solutions) of the differential
equation are analytic functions of $t$; and more precisely,
to determine for a given solution at which values of $t$
it is analytic, and at which values it is not analytic.
Both sufficient conditions for analyticity, and also for
nonanalyticity, at certain values of $t$ are obtained.
It is shown that for some equations there exists a solution
which is $C^\infty$ everywhere, and is analytic at
certain values of $t$ but is not analytic at other values of
$t$. Throughout our analysis, the dynamic properties of the map
$t\to\eta(t)$ play a crucial role.
\vfill

\end{abstract}

\section{Introduction}

We study analyticity properties of solutions of the differential equation
\be{nl}
\dot x(t)=f(t,x(t),x(\eta(t)))
\ee{nl}
with a single variable time-shift $\eta(t)$.
In the broadest setting, we assume the function $f(t,u,v)$ is analytic (that is, holomorphic) in some
region $U$, with $(t,u,v)\in U\subseteq\IR\times\IR^N\times\IR^N$,
and $f:U\to\IR^N$. (Of course $f$ extends to be analytic in a complex neighborhood $\wt U\subseteq\IC\times\IC^N\times\IC^N$
of each point $(t,u,v)\in U$.) We also assume that $\eta(t)$ is a given function
that is analytic in
$t$ in an appropriate region. Our interest is in the analyticity properties
of solutions of equation~\eqref{nl}. In particular, we shall show that
very often periodic solutions, and other globally defined and bounded solutions,
can fail to be analytic at certain points.

In practice the equations we consider will generally be delay-differential equations
with a variable delay $r(t)\ge 0$, that is,
\be{del}
\dot x(t)=f(t,x(t),x(t-r(t))).
\ee{del}
However, we may consider $p$-periodic solutions
$x(t)$ of equation~\eqref{del} when $f(t,u,v)$ and $r(t)$ are $p$-periodic in $t$,
and we note that such solutions are also solutions of equation~\eqref{nl}
when $\eta(t)=t-r(t)+mp$ for any integer $m$. In this fashion, it is natural to consider
equation~\eqref{nl} for which $t-\eta(t)$ takes on both positive and negative
values. As such, we shall study solutions of equation~\eqref{nl} near points
$t_0$ where $\eta(t_0)=t_0$.

State-dependent problems, such as those of the form
\be{state}
\dot x(t)=f(t,x(t),x(t-r)),\qquad r=r(x(t)),
\ee{state}
are certainly of great interest, and provide motivation
for studying the simpler problem~\eqref{del}; see, for example,~\cite{saw} and the references therein. However, we
do not consider such equations in the present paper.

Classic results of one of the authors~\cite{nuss} show that for a broad class of
equations with analytic $f$, but with constant delays, many solutions are analytic.
As a special case, results in~\cite{nuss} show that if $f:\IR^{N(M+1)}\to\IR^N$ is analytic and $r_k\ge 0$
for $1\le k\le M$ are given constants, then any solution $x(t)$ of
\be{1a}
\dot x(t)=f(x(t),x(t-r_1),x(t-r_2),\ldots,x(t-r_M))
\ee{1a}
which exists and is bounded on some interval $(-\infty,t_0]$ must be analytic in $t$.
In fact, these same conclusions also hold when $x(t-r_k)$ is replaced by a distributed delay, namely by
$$
\int_{-r_k}^0x(t+s)\:d\mu_k(s),
$$
where $\mu_k$ is a given signed Borel measure; and this case was the proximate motivation for~\cite{nuss}.
The methods of~\cite{nuss} allow for a straightforward extension to the case of a nonautonomous
periodic system such as
\be{mult}
\dot x(t)=f(t,x(t),x(t-r_1),x(t-r_2),\ldots,x(t-r_M)),
\ee{mult}
where $f(t+p,u_0,u_1,u_2,\ldots,u_M)=f(t,u_0,u_1,u_2,\ldots,u_M)$ holds identically for some $p>0$,
but where again the delays $r_k$ are constant. (We shall study this, and other extensions
of the results of~\cite{nuss}, in a forthcoming paper~\cite{analytic}.)

Until now, however, there have been almost no results for systems with variable delays.
But we mention that analyticity and other properties of certain solutions of the so-called pantograph equation
$$
\dot x(t)=\sum_{k=1}^M\bigg(a_kx(\lam_kt)+b_k\dot x(\lam_kt)\bigg)
$$
and its generalizations, where typically $|\lam_k|\le 1$, have
been studied, among others, by Derfel and Iserles~\cite{derfis},
Iserles~\cite{is93}, and Iserles and Liu~\cite{is94}, \cite{is97a}, \cite{is97}.

We also mention an example shown to us and described in several public lectures by Tibor Krisztin~\cite{tk}, namely
$$
\dot x(t)=f(x(t),x(t-r_1),x(t-r_2),\ldots,x(t-r_M)),\qquad r_k=r_k(x(\cdot)),
$$
where the state-dependent delays $r_k$ are determined implicitly by
$$
\int^t_{t-r_k}a(x(s))\:ds=\rho_k
$$
for $1\le k\le M$. Here $a:\IR\to\IR$ is a given analytic function with $a(u)>0$ for all $u$,
the function $f$ is real-analytic in an appropriate region in $\IR^{M+1}$,
and $\rho_k$ for $1\le k\le M$ are given positive constants. Again, solutions
$x(t)$ which exist and are bounded on some interval $(-\infty,t_0]$ must be analytic.
We show here that this problem in fact can be reduced to one of constant delays
(although we note that Krisztin's original proof of analyticity used a different approach).
First introduce the new time variable
$$
\wt t=\int_{t_0}^t a(x(s))\:ds
$$
and let $y(\wt t)=x(t)$ for $t\le t_0$, equivalently, for $\wt t\le 0$. With $x(t)$
a solution as indicated, then
$$
\dot y(\wt t)=a(y(\wt t))^{-1}f(y(\wt t),y(\wt t-\rho_1),y(\wt t-\rho_2),\ldots,y(\wt t-\rho_M)).
$$
This is an equation with constant delays, and thus $y(\wt t)$
is analytic in $\wt t$ by the classic results in~\cite{nuss}. Reversing the change of variables,
namely letting
$$
t=t_0+\int_0^{\wt t}a(y(s))^{-1}\:ds,
$$
which is an analytic change of variables, shows directly that $x(t)$ is analytic in $t$.

In the present paper, we show that although there are situations where equation~\eqref{nl}
admits nontrivial analytic solutions (see Theorem~2.1), there are also many
robust situations where solutions fail to be analytic (see Theorem~4.2).
The solutions $x(t)$ we consider are often defined for all $t\in\IR$, and in many cases are periodic in $t$,
and thus are $C^\infty$.
In fact very often for such solutions, there can be a coexistence
of points of analyticity and of nonanalyticity; namely, for a given solution $x(t)$,
both the set $\mathcal{A}\subseteq\IR$ of $t$ at which this solution is analytic, and the complement
$\mathcal{N}\subseteq\IR$ of that set,
can be nonempty. It turns out that certain dynamical properties of the mapping $t\to\eta(t)$
on $\IR$ (or on $S^1$ in the case of periodic solutions)
are relevant to the problem of determining the sets of analyticity and of nonanalyticity.

This paper is organized as follows. In Section~2 we consider the (analytic) equation~\eqref{nl}
with initial condition $y(t_0)=y_0$ where $\eta(t_0)=t_0$, and in Theorem~2.1 we show, for the so-called
contractive case $|\dot\eta(t_0)|<1$, that there exists a unique local $C^1$ solution which additionally is analytic.
This result is extended there to the case of a contractive periodic point, that is, $\eta^M(t_0)=t_0$
and $|\dot\eta^M(t_0)|<1$ for some $M>1$.
In Section~3 we study the global mapping properties of the map $\eta$ on the
sets $\mathcal{A}$ and $\mathcal{N}$ of analyticity and nonanalyticity. It is shown in Theorems~3.1 and~3.4
that under mild conditions, and except for a discrete set of exceptional points,
$\eta$ maps $\mathcal{A}$ into $\mathcal{A}$ and $\mathcal{N}$ into $\mathcal{N}$.
In Section~4 we restrict to the scalar linear equation
$$
\dot x(t)=a(t)x(t)+b(t)x(\eta(t))+h(t)
$$
with analytic coefficients $a$, $b$, and $h$, and time-shift $\eta$, in the so-called expansive case
where $\eta(t_0)=t_0$ and $|\dot\eta(t_0)|>1$. For any initial condition $y(t_0)=y_0$ we obtain
a quantity $w_\infty$, depending on $y_0$, such that this initial value problem has a local solution which is analytic
if and only if $w_\infty=0$; these results are found in Theorems~4.2 and~4.5. In Theorem~4.4 we
show that irrespective of the value of $w_\infty$, the initial value problem in fact has an infinite-dimensional
set of local solutions which are $C^\infty$ and not analytic.
Section~5 considers a class of examples of linear equations, given as the integral equation
$$
\kappa x(t)=\int^t_{t-r(t)}x(s)\:ds,\qquad
r(t)=-(\lam-1)\sin t+2\pi m,
$$
where $\kappa$ is an eigenvalue to be determined as part of the solution.
Under appropriate conditions there exists a unique positive $2\pi$-periodic solution $x(t)$ (up to scalar multiple).
It is shown in Theorem~5.2 that under additional condition, $w_\infty\ne 0$ for this solution at
the expansive point $t_0=0$ with $\eta(t)=t+(\lam-1)\sin t$, and thus $x(t)$ is not analytic in any neighborhood of $t=0$.
(This entails rewriting the integral equation as a differential equation.)
In Theorem~5.4 we obtain further conditions under which there additionally exists a contractive point $t_{00}$,
in the neighborhood of which the solution is analytic; thus $x(t)$, which is $C^\infty$ everywhere,
is analytic at some but not all values of $t$, so-called coexistence of analyticity and nonanalyticity.
Finally, in Section~6 we mention several open problems arising from our investigations.

\section{Contractive Periodic Points of {\boldmath $\eta$}: Analyticity}

Our first result is the following, which gives
conditions for the existence of a solution which is analytic
in $t$, at least for a certain range of $t$. One may describe conditions~\eqref{con} and~\eqref{con2}
by saying that the fixed point or periodic point $t_0$ of $\eta$ is {\bf contractive}.
\vv

\noindent {\bf Theorem~2.1.} \bxx
Consider equation~\eqref{nl} where $f:U\to\IR^N$ is analytic in a neighborhood $U\subseteq\IR\times\IR^N\times\IR^N$
of some point $(t_0,x_0,x_0)$, and where $\eta:V\to\IR$ is analytic in a neighborhood $V\subseteq\IR$ of $t_0$.
Assume that $\eta(t_0)=t_0$ and that
\be{con}
|\dot\eta(t_0)|<1.
\ee{con}
Then there exists a unique $C^1$ solution of~\eqref{nl}
with the initial condition $x(t_0)=x_0$ on some interval about $t=t_0$. Moreover, this solution
is analytic in $t$.

The corresponding result holds for periodic points of $\eta$ instead of fixed points.
That is, assume for some $M>1$ that there exist $M$ distinct times $t_n\in\IR$ and $M$ points $x_n\in\IR^N$
for $0\le n\le M-1$, with $(t_n,x_n,x_{n+1})\in U$, with $t_n\in V$,
and with $\eta(t_n)=t_{n+1}$, where we write $t_M=t_0$ and $x_M=x_0$; also assume that
\be{con2}
|\dot\eta^M(t_0)|<1
\ee{con2}
where $\eta^M$ denotes the $M^{\xx th}$ iterate of the function $\eta$.
Then there exists a unique $C^1$ solution of~\eqref{nl}
with the simultaneous initial conditions $x(t_n)=x_n$ on some intervals about $t=t_n$ for each $n$. Moreover, this solution
is analytic in $t$.
\exx
\vv

\noindent {\bf Remark.}
We shall generally let $\eta^n$ denote the $n^{\rm th}$ iterate of the function $\eta$ and $\dot\eta^n(t)=\frac{d}{dt}\eta^n(t)$,
as in the above result. By contrast, $\eta^{(n)}$, as is used for example in Theorem~3.4 below,
denotes the $n^{\rm th}$ derivative of the function~$\eta$.
\vv

\noindent {\bf Proof of Theorem~2.1.}
In a standard fashion we write the differential equation in integrated form
\be{cm}
x(t)=x_0+\int_{t_0}^tf(s,x(s),x(\eta(s)))\:ds
\ee{cm}
and obtain a solution via a contraction mapping argument. In fact, we do this
twice, with two different Banach spaces, namely the space
$$
X=\{x:\overline{D_\delta(t_0)}\to\IC^N\:|\:x(\cdot)\hbox{ is analytic in }D_\delta(t_0)\hbox{ and continuous in }\overline{D_\delta(t_0)}\}
$$
where $D_\delta(t_0)=\{t\in\IC\:|\:|t-t_0|<\delta\}$, and also the more standard space $Y=C([t_0-\delta,t_0+\delta];\IR^N)$.
The supremum norm is taken for each space and the right-hand side of~\eqref{cm} is regarded as a map in the $\eps$-ball
$$
\{x\in X\:|\:\|x-x_0\|\le\eps\}\qquad\mbox{or}\qquad
\{x\in Y\:|\:\|x-x_0\|\le\eps\}.
$$
In particular, the condition~\eqref{con} ensures that if $\delta$ is small enough then
$|\eta(t)-t_0|<|t-t_0|$ when $|t-t_0|\le\del$; and for appropriately chosen $\del$ and $\eps$
the relevant mapping is well-defined and a contraction
from the above balls into themselves. It follows that the fixed point in $X$ gives an analytic solution,
and this solution is unique among all elements of $Y$, as desired. We omit the details.

For the case when $\eta$ has a periodic point, let
$$
y_n(t)=x(\eta^n(t_0+t)),\qquad
\mbox{for }0\le n\le M-1.
$$
Then equation~\eqref{nl} can be written as the system
\be{nl2}
\dot y_n(t)=\left\{\begin{array}{ll}
f_n(t,y_n(t),y_{n+1}(t)), &\bb \quad\mbox{for } 0\le n\le M-2,\\
\\
f_n(t,y_n(t),y_0(\eta^M(t_0+t)-t_0)), &\bb \quad\mbox{for } n=M-1,
\end{array}\right.
\ee{nl2}
now in a neighborhood of $t=0$, with nonlinearities
$$
f_n(t,u,v)=\dot\eta^n(t_0+t)f(\eta^n(t_0+t),u,v).
$$
This reduces the problem
to the fixed-point case considered above, now with the fixed point $\wt\eta(0)=0$ and $|\dot{\wt\eta}(0)|<1$
where $\wt\eta(t)=\eta^M(t_0+t)-t_0$.~\fp
\vv

If $|\dot\eta(t_0)|>1$ in place of condition~\eqref{con} in Theorem~2.1, then we have
a fixed point of $\eta(t)$ which we may term {\bf expansive}. As will be shown in Section~4, a linear
equation with an expansive fixed point generally does not possess an analytic solution
through that point (although it has many $C^\infty$ solutions there).
In fact, we shall show that typically
(in a sense to be made precise) no analytic solution exists about an expansive fixed point,
although there are exceptional equations which do possess one.

\section{Mapping Properties of {\boldmath $\eta$}}

Here we examine the role of the mapping $t\to\eta(t)$ in determining regions of analyticity and nonanalyticity
for solutions of~\eqref{nl}. Generally, if $x(t)$ is such a solution which is defined for all $t\in\IR$,
we define sets
\be{an}
\begin{array}{lcl}
\mathcal{A} &\bb = &\bb
\{t_0\in\IR\:|\:x(t)\hbox{ is analytic for }|t-t_0|<\del,\hbox{ for some }\del\},\\
\\
\mathcal{N} &\bb = &\bb
\IR\setminus\mathcal{A}.
\end{array}
\ee{an}
Clearly, $\mathcal{A}$ is an open set while $\mathcal{N}$ is closed.
As Theorems~3.1 and~3.4 below show, the dynamic properties of the map $t\to\eta(t)$ are key, and
in particular, the sets $\mathcal{A}$ and $\mathcal{N}$ enjoy certain mapping properties with respect to $\eta$.

The case of a periodic solution of equation~\eqref{nl} is certainly of great interest. Here one has
\be{xper}
x(t+p)=x(t)
\ee{xper}
for all $t\in\IR$, for some $p>0$, along with the conditions
\be{per}
f(t+p,u,v)=f(t,u,v),\qquad
\eta(t+p)=\eta(t)+p,
\ee{per}
for all $t\in\IR$ on the differential equation.
(Even in the case that $f(t,u,v)$ is linear
in $u$ and $v$, such solutions can arise naturally as Floquet solutions.)
Here we may regard $\eta$ as a map $\eta:S^1\to S^1$ of the
circle $S^1=\IR/p\ZZ$ onto itself.

The following elementary result for analytic differential equations shows that if
a solution $x(t)$ is analytic at $\eta(t_0)$ for some $t_0$ in its domain, then it
is analytic at $t_0$.
\vv

\noindent {\bf Theorem~3.1.} \bxx
Assume that $\eta:V_0\to\IR$ is analytic in some neighborhood $V_0\subseteq\IR$ of a point $t_0$,
and denote $t_1=\eta(t_0)$. Also assume for some neighborhood $V_1\subseteq\IR$ of $t_1$
that $x:V_0\cup V_1\to\IR^N$ is a continuous function. In addition, assume that $f:U\to\IR^N$
is analytic in some neighborhood $U\subseteq\IR\times\IR^N\times\IR^N$ containing the point $(t_0,x(t_0),x(t_1))$, and that $x(t)$
is $C^1$ and satisfies the differential equation~\eqref{nl} for all $t\in\IR$ near $t_0$.
Then if $x(t)$ is analytic in a neighborhood of $t_1$, it is analytic in a neighborhood of $t_0$.
\exx
\vv

\noindent {\bf Proof.}
Assuming that $x(t)$ is analytic in $t$ in a neighborhood of $t=t_1$, we have that
$x(\eta(t))$ is analytic in $t$ in a neighborhood of $t=t_0$. Upon regarding $x(\eta(t))$
as a known function and equation~\eqref{nl} as an ordinary differential equation for $x(t)$,
it is then immediate
that the solution $x(t)$ of this equation is analytic for $t$
in a neighborhood of $t_0$, as desired.~\fp
\vv

Theorem~3.1 can be used to provide the following global result, involving iterates of $\eta$.
\vv

\noindent {\bf Corollary~3.2.} \bxx
Consider equation~\eqref{nl} where $f:U\to\IR^N$ is analytic in a region $U\subseteq\IR\times\IR^N\times\IR^N$
and where $\eta:\IR\to\IR$ also is analytic.
Assume that $x(t)$ is a solution of~\eqref{nl} for all $t\in\IR$, satisfying
$(t,x(t),x(\eta(t)))\in U$ for all $t$.
Also assume there exists a point $t_0\in\IR$ such that $\eta^M(t_0)=t_0$ and $|\dot\eta^M(t_0)|<1$ for some $M\ge 1$, let
$t_k=\eta^k(t_0)$ for $1\le k\le M-1$, and let
$$
\mathcal{B}=\bigcup_{k=0}^{M-1}\mathcal{B}_k,\qquad
\mathcal{B}_k=\{t\in\IR\:|\:\lim_{n\to\infty}\eta^{Mn}(t)=t_k\},
$$
which is the basin of attraction of the orbit of the periodic point $t_0$ for the map $\eta$.
Then $\mathcal{B}\subseteq\mathcal{A}$, that is, $x(t)$ is analytic at every point of $\mathcal{B}$.

Suppose additionally for some $p>0$ that the periodicity conditions~\eqref{per}
hold, and that $x(t)$ is $p$-periodic, namely~\eqref{xper} holds, for all $t\in\IR$.
Also relax the condition $\eta^M(t_0)=t_0$ to assume simply $\eta^M(t_0)=t_0 \:({\rm mod}\:p)$,
but keep the condition $|\dot\eta^M(t_0)|<1$. Similarly let $\mathcal{B}\subseteq S^1=\IR/p\ZZ$ denote
the basin of attraction of the orbit of $t_0$ with $\eta:S^1\to S^1$ considered as a map on the circle. Then
again, $\mathcal{B}\subseteq\mathcal{A}$, that is, $x(t)$ is analytic at every point of $\mathcal{B}$.
\exx
\vv

\noindent {\bf Proof.}
By Theorem~2.1 the solution $x(t)$ is analytic in some neighborhood $V\subseteq\IR$ of $\{t_0,t_1,\ldots,t_{M-1}\}$.
Let $\tau\in\mathcal{B}$, say $\tau\in\mathcal{B}_k$ where $0\le k\le M-1$.
Then $\eta^{Mn}(\tau)\in V$ for some integer $n\ge 0$, thus $x(t)$ is analytic
in a neighborhood of $\eta^{Mn}(\tau)$.
It follows by Theorem~3.1 that $x(t)$ is analytic in some neighborhood of $\tau$,
and so $\tau\in\mathcal{A}$, and thus $\mathcal{B}\subseteq\mathcal{A}$, as desired.

The time-periodic case is proved with minor modifications.~\fp
\vv

Theorem~3.4 below provides a partial converse to Theorem~3.1, but requires additional conditions on $f$
(for example, $f(t,u,v)$ must have some nontrivial dependence on $v$).
Even in the scalar case $N=1$ these conditions are rather intricate; we describe them here.
We first define polynomials $P_n$, for $n\ge 0$, inductively by
\be{pol}
\begin{array}{l}
P_0(\zet_{0,0})=\zet_{0,0},\\
\\
\ds{P_n(\{\zet_{i,j}\}_{i+j\le n})
=\sum_{k+m\le n-1}\bigg(D_{\zet_{k,m}}P_{n-1}(\{\zet_{i,j}\}_{i+j\le n-1})\bigg)
\bigg(\zet_{k+1,m}+\zet_{0,0}\zet_{k,m+1}\bigg).}
\end{array}
\ee{pol}
The scalar variables $\zet_{i,j}$, for nonnegative integers $i$ and $j$ in the range $i+j\le n$, are the arguments
of $P_n$, and thus $P_n$ is a polynomial in $\frac12 (n+1)(n+2)$ variables.
The interpretation of the formula~\eqref{pol} is as follows. If $x(t)$ is a solution of
a scalar ordinary differential equation $\dot x(t)=f(t,x(t))$, then upon repeated differentiation with
respect to $t$ and substitution of the differential equation into the formula obtained,
we arrive at a formula for the derivative $x^{(n+1)}(t)$ which is a polynomial in
the derivatives $D_t^iD_x^jf(t,x(t))$. The polynomial so obtained
is simply $P_n$ with the substitution of $\zet_{i,j}=D_t^iD_x^jf(t,x(t))$
for its arguments. In particular,
$$
P_1(\{\zet_{i,j}\}_{i+j\le 1})=\zet_{1,0}+\zet_{0,0}\zet_{0,1}
$$
corresponds to the formula $\ddot x(t)=D_tf(t,x(t))+f(t,x(t))D_xf(t,x(t))$, and
$$
P_2(\{\zet_{i,j}\}_{i+j\le 2})=
\zet_{2,0}+\zet_{0,0}\zet_{1,1}
+(\zet_{1,0}+\zet_{0,0}\zet_{0,1})\zet_{0,1}
+\zet_{0,0}(\zet_{1,1}+\zet_{0,0}\zet_{0,2})
$$
corresponds to the analogous formula for $x^{(3)}(t)$ (which we omit).

With this, we have the following lemma in the case of a scalar equation.
\vv

\noindent {\bf Lemma~3.3.} \bxx
Assume that $f:U\to\IR$ is analytic in some neighborhood $U\subseteq\IR\times\IR\times\IR$
of a point $(t_0,x_0,v_0)$, and define the functions
$$
Q_n(v)=D_v\bigg(P_n(\{D_t^iD_x^jf(t_0,x_0,v)\}_{i+j\le n})\bigg)
$$
for $v$ near $v_0$.
Then there exists an analytic
function $x(t)$, for $t$ near $t_0$ and with $x(t_0)=x_0$, such that
\be{g}
g(t,v)=\dot x(t)-f(t,x(t),v)
\ee{g}
is the zero function in a neighborhood of $(t_0,v_0)$, if and only if each
of the functions $Q_n(v)$ is identically zero in a neighborhood of $v_0$.
\exx
\vv

\noindent {\bf Proof.}
Suppose first there exists an analytic function $x(t)$ as in the statement
of the lemma. Then for every $v$ near $v_0$ this function satisfies the equation $\dot x(t)=f(t,x(t),v)$
with the initial value $x(t_0)=x_0$, and so
\be{xp}
x^{(n+1)}(t_0)=P_n(\{D_t^iD_x^jf(t_0,x_0,v)\}_{i+j\le n})
\ee{xp}
for each $n\ge 0$, from the remarks following the definition of $P_n$ above. The left-hand
side of~\eqref{xp} is independent of $v$, and so the derivative of the right-hand side, namely
$Q_n(v)$, vanishes identically in a neighborhood of $v_0$, as claimed.

Now suppose that $Q_n(v)=0$ identically in a neighborhood of $v_0$, for each $n\ge 0$.
For $v$ in a neighborhood of $v_0$
let $y(t,v)$ be the solution of the ordinary differential equation $\dot y(t,v)=f(t,y(t,v),v)$
with the initial value $y(t_0,v)=x_0$. Also let $x(t)=y(t,v_0)$ and define $g(t,v)$ as in~\eqref{g}.
Then again from the definition of $P_n$
$$
D_t^{n+1}y(t_0,v)=P_n(\{D_t^iD_x^jf(t_0,x_0,v)\}_{i+j\le n})
$$
for every $n\ge 0$ and so
$$
D_vD_t^{n+1}y(t_0,v)=D_v\bigg(P_n(\{D_t^iD_x^jf(t_0,x_0,v)\}_{i+j\le n})\bigg)=Q_n(v)=0
$$
holds identically for $v$ in some neighborhood of $v_0$. Therefore, the quantity
$D_t^{n+1}y(t_0,v)$ is constant in $v$, and so it equals its value at $v=v_0$, namely
$$
D_t^{n+1}y(t_0,v)=x^{(n+1)}(t_0),
$$
for every $n\ge 0$.
As both $y(t,v)$ and $x(t)$ are analytic in $t$, it follows they are equal as they have the same Taylor series about $t_0$.
Therefore, for any $(t,v)$ near $(t_0,v_0)$ one has
$$
g(t,v)=\dot x(t)-f(t,x(t),v)=D_ty(t,v)-f(t,y(t,v),v)=0
$$
and so $g$ is the zero function, as claimed.~\fp
\vv

\noindent {\bf Theorem~3.4.} \bxx
Assume all the conditions in the statement of Theorem~3.1, except for the final sentence.
Also assume
that $x(t)$ is $C^\infty$ in a neighborhood of $t_1$ and that one of the following two conditions holds:
\begin{itemize}
  \item[{(1)}] $N=1$ (a scalar equation), and there exists $n\ge 0$ such that the function $Q_n(v)$, with $x_0=x(t_0)$,
is not identically zero in a neighborhood of $v_0=x(\eta(t_0))$; or
  \item[{(2)}] $f(t,u,v)=f_0(t,u)+B(t)v$ in a neighborhood of $(t_0,x(t_0),x(\eta(t_0)))$,
and the function $\det B(t)$ does not vanish identically near $t_0$.
\end{itemize}
Finally assume that there exists an integer $m\ge 0$ such that
$$
\eta^{(k)}(t_0)=0\mbox{ for }1\le k\le 2m,\qquad\eta^{(2m+1)}(t_0)\ne 0,
$$
that is, either $\dot\eta(t_0)\ne 0$ or else $\dot\eta(t)$ has a zero of even order at $t_0$.
Then if $x(t)$ is analytic in a neighborhood of $t_0$, it is analytic in a neighborhood of $t_1$.
\exx
\vv

\noindent {\bf Proof.}
Assume that $x(t)$ is analytic for $t$ in a neighborhood of $t_0$. First suppose that
condition~(1) holds, and define
$$
g(t,v)=\dot x(t)-f(t,x(t),v),\qquad
\nu(t)=x(\eta(t)).
$$
Then $g(t,v)$ is analytic for $(t,v)$ in a neighborhood of $(t_0,x(t_1))$, and
$g(t,\nu(t))=0$ identically for $t$ in a neighborhood of $t_0$. Lemma~3.3 implies that $g(t,v)$
is not the zero function, and this fact together with a straightforward application
of Newton's polygon (see, for example,~\cite[Section~2.8]{chha} or~\cite[Theorem~1.4, Section~1.7]{eb})
implies that $\nu(t)$ is given by a fractional power series
\be{fs}
\nu(t)=\sum_{j=0}^\infty \nu_j(t-t_0)^{j/q},
\ee{fs}
for some integer $q\ge 1$, which converges for $t$ in some neighborhood of $t_0$. However, $\nu(t)$
is a $C^\infty$ function of $t$ near $t_0$, since $x(t)$ is $C^\infty$ near $t_1$. Necessarily then, the only
nonzero terms in the power series~\eqref{fs} are those for which $j=iq$ is a multiple of $q$.
That is,
$$
\nu(t)=\sum_{i=0}^\infty \nu_{iq}(t-t_0)^i,
$$
and we conclude that $\nu(t)$ is analytic for $t$ in a neighborhood of $t_0$.

Next, we have that $\eta(t)=\eta(t_0)+((t-t_0)\theta(t))^{2m+1}=t_1+((t-t_0)\theta(t))^{2m+1}$ for some function
$\theta(t)$ which is analytic in a neighborhood of $t_0$ with $\theta(t_0)\ne 0$. As the function
$\nu(t)=x(t_1+((t-t_0)\theta(t))^{2m+1})$ is analytic near $t=t_0$, it follows via
the analytic change of variables $\wt t=(t-t_0)\theta(t)$
that the function $\mu(\wt t)=x(t_1+\wt t^{2m+1})$ is analytic in $\wt t$ in a neighborhood of $\wt t=0$. Writing
$$
\mu(\wt t)=\sum_{j=0}^\infty \mu_j\wt t^j,
$$
which is a convergent series, we conclude that
\be{fs2}
x(t)=\mu((t-t_1)^{1/(2m+1)})=\sum_{j=0}^\infty \mu_j(t-t_1)^{j/(2m+1)}
\ee{fs2}
for $t$ in a neighborhood of $t_1$. Again, as $x(t)$ is $C^\infty$, the only nonzero terms in~\eqref{fs2}
are those for which $j=i(2m+1)$ is a multiple of $2m+1$. It follows that $x(t)$ is analytic for $t$ in a neighborhood
of $t_1$, as desired. (In particular, the oddness of the denominator $2m+1$ implies the formula~\eqref{fs2} is
valid both for $t>t_1$ and $t<t_1$.)

In case condition~(2) is assumed in place of~(1), the proof is similar. In particular, the function $\nu(t)$
has the form
$$
\nu(t)=B(t)^{-1}(\dot x(t)-f_0(t,x(t)))
$$
for $t\ne t_0$, with a possible pole at $t=t_0$, and is thus a meromorphic function of $t$ in a neighborhood of $t_0$.
However, as $\nu(t)$ is continuous at $t_0$, we have that
in fact $\nu(t)$ is analytic in $t$. The remainder of the proof follows as above.~\fp
\vv

\noindent {\bf Remark.}
For the case $N=1$ in Theorem~3.4, condition~(1) is sharp in the sense that it gives necessary and
sufficient conditions (as per Lemma~3.3) for the existence of $x(t)$ so that the function $g(t,v)$ is identically zero.
As a practical matter, the condition for $n=0$, namely that $Q_0(v)=D_vf(t_0,x(t_0),v)$ does not vanish identically, or equivalently
that $f(t_0,x(t_0),v)$ has nontrivial dependence on $v$,
should suffice for many situations. And if $Q_0(v)$ does vanish identically, one could next require that
$$
Q_1(v)=D_v\bigg(D_tf(t_0,x(t_0),v)+f(t_0,x(t_0),v)D_xf(t_0,x(t_0),v)\bigg)
$$
does not vanish identically.

More generally, noting that $Q_n(v)=Q_n(v;t_0,x_0)$ depends on the choice of $(t_0,x_0)$,
one could (working in $U=\IR\times\IR^N\times\IR^N$ for simplicity) define sets
$$
\begin{array}{lcl}
S_n &\bb = &\bb \{(t_0,x_0)\in\IR\times\IR^N\:|\:Q_k(v;t_0,x_0)=0\mbox{ for every }v\in\IR^N\mbox{ and }0\le k\le n\},\\
\\
S_\infty &\bb = &\bb \ds{\bigcap_{n=0}^\infty S_n,\qquad
T=\{t_0\in\IR\:|\:\mbox{there exists }x_0\in\IR^N\mbox{ with }(t_0,x_0)\in S_\infty\}.}
\end{array}
$$
One could then require that the set $T\subseteq\IR$ be a totally disconnected (either discrete or at least zero-dimensional) set,
and observe that if $t_0\not\in T$ then condition~(1) of Theorem~3.4 applies at every $x_0\in\IR^N$. One would
expect the theory of semianalytic and subanalytic sets to play a role here.
\vv

\noindent {\bf Remark.}
It is not clear what are the appropriate conditions in place of~(1) or~(2) in Theorem~3.4
for a general nonlinearity $f(t,u,v)$ with $N>1$. If $N=1$ then the necessary and sufficient
condition on an analytic function $g(t,v)$ in order to conclude that any $C^\infty$ solution
$\nu(t)$ to $g(t,\nu(t))=0$ is analytic, as in the above proof, is that the function $g$
does not vanish identically. But for $N>1$ this condition is not sufficient; for example,
if $g(t,v)=Bv$ where $B$ is a constant matrix for which $\det B=0$ but $B\ne 0$, then there exist solutions
$\nu(t)$ to $B\nu(t)=0$ which are $C^\infty$ but not analytic. Necessary and sufficient
conditions on $g$ in order to conclude that $\nu(t)$ is analytic have in fact been given by Neelon~\cite{neelon},
and involve an iterated Jacobian ideal in the ring of germs of analytic functions at $(t_0,v_0)$.
Neelon's result uses in a crucial way the so-called Artin Approximation Theorem~\cite{artin}, which states
that if $\nu(t)$ is a formal power series solution of $g(t,\nu(t))=0$, where $g(t,v)$ is analytic
in $(t,v)$, then for every $n\ge 0$ there exists a true solution $\wt\nu(t)$ which is a convergent
power series which agrees with $\nu(t)$ up to order $n$ in powers of $t$.
\vv

One immediate consequence of Theorems~3.1 and~3.4 is the following result.
Note that in any case, the solution $x(t)$ in this result is everywhere $C^\infty$ by virtue of its
existing for all time.
\vv

\noindent {\bf Corollary~3.5.} \bxx
Consider equation~\eqref{nl} where $f:U\to\IR^N$ is analytic in a region $U\subseteq\IR\times\IR^N\times\IR^N$
and where $\eta:\IR\to\IR$ also is analytic.
Assume that $x(t)$ is a solution of~\eqref{nl} for all $t\in\IR$, satisfying
$(t,x(t),x(\eta(t)))\in U$ for all $t$.
Also assume that at every $t_0\in\IR$ either condition~(1) of Theorem~3.4 holds (where $n$ may depend on $t_0$)
or else that condition~(2) of Theorem~3.4 holds. Recall the sets $\mathcal{A}$ and $\mathcal{N}$
of analyticity and nonanalyticity, respectively, in~\eqref{an},
and also define the set
$$
\begin{array}{lcl}
\mathcal{M} &\bb = &\bb
\{t_0\in\IR\:|\:\mbox{there exists an integer }m>0\mbox{ such that}\\
\\
&\bb &\bb
\eta^{(k)}(t_0)=0\mbox{ for }1\le k\le 2m-1\mbox{ but }\eta^{(2m)}(t_0)\ne 0\}.
\end{array}
$$
Then
$$
\eta(\mathcal{A}\setminus\mathcal{M})\subseteq\mathcal{A},\qquad
\eta(\mathcal{N})\subseteq\mathcal{N},
$$
both hold. In particular, if for all $t\in\IR$ one has
\be{mon}
\dot\eta(t)\ge 0,\qquad\eta(\pm\infty)=\pm\infty,
\ee{mon}
then $\mathcal{M}=\emptyset$ and
\be{3a}
\eta(\mathcal{A})=\mathcal{A},\qquad
\eta(\mathcal{N})=\mathcal{N}.
\ee{3a}
In this case $\eta(t_0)\in\mathcal{A}$ if and only if $t_0\in\mathcal{A}$, and similarly for $\mathcal{N}$,
for every $t_0\in\IR$.
\exx
\vv

\noindent {\bf Proof.}
The results follow directly from Theorems~3.1 and~3.4.~\fp
\vv

If, in the periodic case~\eqref{per} of the above result, the monotonicity condition~\eqref{mon} holds
for every $t\in\IR$, then the map $\eta$ can be regarded as a homeomorphism of the circle
$S^1=\IR/p\ZZ$ onto itself.
(The analyticity of $\eta$ implies it is one-to-one.)
Thus the rotation number $\omega$ of $\eta$ is defined, namely
\be{rot}
\omega=\lim_{n\to\pm\infty}\frac{\eta^n(t_0)}{np},
\ee{rot}
where the limit exists and is independent of the sign $\pm$ and the choice of $t_0$.
If $\omega$ is rational
then there exists a periodic point (modulo $p$) of the map $\eta$.
If $\omega$ is an integer then necessarily there exists a fixed point (modulo $p$) of $\eta$,
and for generic such $\eta$ there exist both contractive and expansive fixed points.
In case $\omega$ is irrational then (following Yoccoz~\cite{yoccoz}) $\eta$ is topologically conjugate to the rigid rotation
$\wt\eta(t)=t+\omega p$ and every orbit $\{\eta^n(t)\}_{n\in\ZZ}$ is dense in the circle.
From this observation we have the following result.
\vv

\noindent {\bf Corollary~3.6.} \bxx
Assume the conditions in the statement of Corollary~3.5, including the monotonicity condition~\eqref{mon} for all $t\in\IR$.
Also assume for some $p>0$ that the periodicity conditions~\eqref{per} hold,
and suppose that $x(t)$ is a solution of~\eqref{nl} of period
$p$, that is, satisfying~\eqref{xper} for all $t\in\IR$.
Finally assume that the rotation number $\omega$ in~\eqref{rot} is irrational.

Then either $x(t)$ is analytic for every $t\in\IR$, or it is not analytic at any $t\in\IR$. In any case, it
is $C^\infty$ at every $t\in\IR$.
\exx
\vv

\noindent {\bf Proof.}
As noted earlier the solution $x(t)$ is $C^\infty$ everywhere. Let $\mathcal{A}$ and $\mathcal{N}$
be as in~\eqref{an}, and recall that these sets
are open and closed, respectively. These sets are also $p$-periodic, namely that $t_0\in\mathcal{A}$ if
and only if $t_0+p\in\mathcal{A}$ and similarly for $\mathcal{N}$. Suppose that $\mathcal{N}\ne\emptyset$ and take
any point $t_0\in\mathcal{N}$. Then by~\eqref{3a} in Corollary~3.5 and by periodicity, we have that $\eta^n(t_0)+mp\in\mathcal{N}$
for every $m,n\in\ZZ$.

The fact that $\omega$ is irrational, together with the analyticity of $\eta$, implies by a result of Yoccoz~\cite{yoccoz}
that $\eta$ is topologically conjugate to the rigid rotation $t\to t+\omega p$.
(This result is classical~\cite{denjoy} --- see also~\cite{hale} --- if $\eta$ is a diffeomorphism,
that is, if $\dot\eta(t)>0$ for every $t\in\IR$. However,
this conclusion can fail even for $C^\infty$ homeomorphisms of $S^1$ which are not diffeomorphisms; see Hall~\cite{hall}.)
Thus the set of points $\eta^n(t_0)+mp$ is dense in $\IR$.
Therefore $\mathcal{N}$ is dense in $\IR$ and so $\mathcal{N}=\IR$ as $\mathcal{N}$ is closed.

We thus conclude that either $\mathcal{N}=\emptyset$ and $\mathcal{A}=\IR$, or else $\mathcal{N}=\IR$ and $\mathcal{A}=\emptyset$,
which is as claimed.~\fp
\vv

\noindent {\bf Remark.}
In~\cite{yoc2} Yoccoz identifies a class $\mathcal{H}\subseteq\IR$ of irrational numbers
with the property that any analytic diffeomorphism $\eta$ of the circle $S^1$
with rotation number $\ome\in\mathcal{H}$ is analytically conjugate to a rigid
rotation. That is, there exists an analytic diffeomorphism $\sig:S^1\to S^1$ such that
$\sig^{-1}(\eta(\sig(t)))=t+\ome p$.
The characterization of $\mathcal{H}$ is sharp, namely that for each $\ome\not\in\mathcal{H}$ there exists
an analytic diffeomorphism $\eta$ with rotation number $\ome$ which is not
analytically conjugate to a rigid rotation.
Further, the set $\mathcal{H}$ properly contains all the Diophantine numbers, namely numbers $\ome$ for
which there exist $K>0$ and $\delta>0$ such that
$|n\ome-m|\ge Kn^{-1-\delta}$ for all integers $m,n\in\ZZ$ with $n\ne 0$.
It follows from this that in the setting of Corollary~3.6, if $\ome\in\mathcal{H}$, then $x(t)$ is analytic for every $t\in\IR$.
Indeed, this is easily shown via the time transformation $t=\sig(\wt t)$ and letting $y(\wt t)=x(\sig(\wt t))$. Then $y(\wt t)$
is a $p$-periodic solution of the analytic equation
$$
\dot y(\wt t)=\dot\sig(\wt t)f(\sig(\wt t),y(\wt t),y(\wt t+\ome p))
$$
with a constant time-shift $\ome p$. (In case $\ome p>0$ one may reverse time $\wt t\to -\wt t$ and transform
this into a delay.) As noted, by a straightforward extension (see~\cite{analytic}) of the results of~\cite{nuss},
the solution $y(\wt t)$, and therefore also $x(t)$, are everywhere analytic.

Whether or not the set $\mathcal{H}$ is sharp with respect to this property of analyticity of periodic
solutions is not clear. Namely, it is unclear whether or not there exists an irrational number $\ome$,
with $\ome\not\in\mathcal{H}$, such that every solution $x(t)$ in the setting of Corollary~3.6 with this $\ome$
must be everywhere analytic.

It would be of great interest to find an example of an analytic equation with a solution as in Corollary~3.6
which is everywhere $C^\infty$ but nowhere analytic.

\section{Expansive Fixed Points of {\boldmath $\eta$}: Nonanalyticity (Usually)}

We next consider the situation in which~\eqref{con} fails; more specifically, we
assume the case of an expansive fixed point, namely that $\eta(t_0)=t_0$ and
\be{4aa}
|\dot\eta(t_0)|>1.
\ee{4aa}
We shall not consider the case of expansive periodic points.
We shall also restrict ourselves to the case of a scalar ($N=1$) linear inhomogeneous system,
which we write as
\be{slin}
\dot x(t)=a(t)x(t)+b(t)x(\eta(t))+h(t),
\ee{slin}
although we expect that analogs of our theorems should hold for expansive periodic points, for systems, and for nonlinear equations.
Thus we assume that $a(t)$, $b(t)$, $h(t)$, and $\eta(t)$ are analytic in $t$ in a neighborhood of $t=t_0$,
with $\eta(t_0)=t_0$ and~\eqref{4aa} holding.

A simple example which provides some insight into the general situation is the equation
\be{simple}
\dot x(t)=a_0x(t)+b_0x(\lambda t).
\ee{simple}
Here we take $a_0,b_0\in\IR$ with $b_0\ne 0$ and $\eta(t)=\lambda t$ with $|\lambda|>1$.
We ask if, for a given $x_0\in\IR$, equation~\eqref{simple} possesses a solution $x(t)$
with $x(0)=x_0$ which is analytic in a neighborhood of $t=0$. Such a solution would take the form
\be{ts}
x(t)=\sum_{n=0}^\infty x_nt^n
\ee{ts}
and one easily checks that the coefficients $x_n$ for $n\ge 1$ are uniquely determined
and given by the formula
$$
x_n=\bigg(\frac{\lambda^{n(n-1)/2}b_0^n}{n!}\bigg)w_n,\qquad
w_n=x_0\prod_{k=0}^{n-1}\bigg(1+\frac{a_0}{\lambda^kb_0}\bigg).
$$
Note that the limit
$$
w_\infty=\lim_{n\to\infty}w_n=x_0\prod_{k=0}^\infty\bigg(1+\frac{a_0}{\lambda^kb_0}\bigg)
$$
exists and is finite. Certainly, if $a_0+\lambda^kb_0=0$ for some $k\ge 0$ then $x_n=w_n=0$ for every $n>k$,
and equation~\eqref{simple}
with $x(0)=x_0$ possesses an analytic solution which is in fact a polynomial. (Of course $w_\infty=0$
in this case.) On the other hand, if $a_0+\lambda^kb_0\ne 0$ for every $k\ge 0$ and if $x_0\ne 0$,
then $w_\infty\ne 0$. In this case
$$
\lim_{n\to\infty}|x_n|^{1/n}=\lim_{n\to\infty}\bigg|\frac{\lambda^{(n-1)/2}b_0}{(n!)^{1/n}}\bigg||w_n|^{1/n}=\infty,
$$
and so the Taylor series~\eqref{ts} has zero radius of convergence and no analytic solution exists. Note that Stirling's formula
$$
n!\sim n^{n+1/2}(2\pi)^{1/2}e^{-n}
$$
is useful in obtaining the above limit. Thus, roughly speaking, for ``most'' but not for every choice of $a_0$ and $b_0$,
equation~\eqref{simple} has no nontrivial analytic solution in a neighborhood of $t=0$.

For the general equation~\eqref{slin} with a given initial condition $x(t_0)=x_0$,
in Theorem~4.2 we shall define a quantity $w_\infty$ analogous
to the one above and show that a necessary condition for an analytic solution in a neighborhood of $t_0$
to exist is that $w_\infty=0$.
In Theorem~4.5 we show this is also sufficient; if $w_\infty=0$ then there does exist such an analytic solution.
We shall also show, in Theorem~4.4, that in any case (either if $w_\infty=0$ or if $w_\infty\ne 0$) there exists
an infinite-dimensional set of solutions with $x(t_0)=x_0$ which are $C^\infty$ but not analytic.

In obtaining these results
it will be useful to transform $\eta(t)$ to the linear function $t_0+\lambda(t-t_0)$,
where $\lam=\dot\eta(t_0)$, in a neighborhood of
$t=t_0$. To this end, the following one-dimensional version of the Hartman--Grobman Theorem
provides such a transformation via a local analytic conjugacy $\sigma(t)$.
\vv

\noindent {\bf Lemma~4.1.} \bxx
Let $\eta:V\to\IR$ be analytic in some neighborhood $V\subseteq\IR$ of
a point $t_0$. Assume that $\eta(t_0)=t_0$ and that $|\lambda|\ne 0,1$ where
$\lambda=\dot\eta(t_0)$.
Then there exists a function $\sigma(t)$, analytic in a neighborhood of $t=t_0$, such that
\be{con-}
\sigma(t_0)=t_0,\qquad\dot\sigma(t_0)=1,\qquad\sigma^{-1}(\eta(\sigma(t)))=t_0+\lambda(t-t_0),
\ee{con-}
holds identically.
\exx
\vv

A proof of Lemma~4.1 may be found in~\cite[Chapter~II, Theorem~2.1]{cargam},
at least for $0<|\lam|<1$; the case $|\lam|>1$ is handled by simply considering
$\eta^{-1}$ in place of $\eta$. We note that in fact this
result was originally given by G.~Koenigs over a century ago~\cite{gk}.
We remark that the function $\sigma$ is unique, not only among analytic
functions, but among all local $C^1$ diffeomorphisms satisfying~\eqref{con-}.
In general $\sigma(t)$ need not be defined outside
a neighborhood of $t=t_0$, although if $\eta(t)$ is an entire function of $t$ and $|\lambda|>1$, then
$\sig(t)$ is also entire. One sees this
by applying $\sig$ to both sides of the final equation in~\eqref{con-}, and then repeatedly iterating.

Making the change of variables
\be{cov}
t=\sigma(t_0+\wt t)
\ee{cov}
in equation~\eqref{slin}, and denoting $y(\wt t)=x(\sigma(t_0+\wt t))$
for $\wt t$ near zero, yields the pantograph-type equation
\be{yeqq}
\dot y(\wt t)=\alpha(\wt t)y(\wt t)+\beta(\wt t)y(\lambda\wt t)+\gamma(\wt t)
\ee{yeqq}
with coefficients
$$
\alpha(\wt t)=\dot\sigma(t_0+\wt t)a(\sigma(t_0+\wt t)),\qquad
\beta(\wt t)=\dot\sigma(t_0+\wt t)b(\sigma(t_0+\wt t)),\qquad
\gamma(\wt t)=\dot\sigma(t_0+\wt t)h(\sigma(t_0+\wt t)),
$$
which are analytic in a neighborhood of $\wt t=0$.
It is certainly the case that $x(t)$ is analytic for $t$ in a neighborhood of $t_0$
if and only if $y(\wt t)$ is analytic for $\wt t$ in a neighborhood of zero.
In what follows we shall use the Taylor series
\be{abg}
\alpha(\wt t)=\sum_{n=0}^\infty \alpha_n\wt t^n,\qquad
\beta(\wt t)=\sum_{n=0}^\infty \beta_n\wt t^n,\qquad
\gamma(\wt t)=\sum_{n=0}^\infty \gamma_n\wt t^n,
\ee{abg}
of the coefficient functions. We note the bounds
\be{the2}
|\alpha_k|,|\beta_k|,|\gamma_k|\le C\mu^k,
\ee{the2}
for some $C>0$ and $\mu>0$ and every $k$, due to the positive
radius of convergence of these functions.

The next result provides a sufficient condition for the absence of an analytic solution
to an initial value problem for equation~\eqref{yeqq}. Observe that the condition $\beta(0)=\beta_0\ne 0$
in this result is equivalent to $b(t_0)\ne 0$ in~\eqref{slin}.
\vv

\noindent {\bf Theorem~4.2.} \bxx
Consider equation~\eqref{yeqq} with analytic coefficients $\alpha,\beta,\gamma:V\to\IR$
in some neighborhood $V\subseteq\IR$ of $\wt t=0$,
and with $|\lambda|>1$. Let $y_0\in\IR$ be given
and define $y_n\in\IR$, for $n\ge 1$, to be the coefficient of $\wt t^n$ in the formal power series for
a solution $y(\wt t)$ of equation~\eqref{yeqq} with initial condition $y(0)=y_0$; that is,
\be{yr}
(n+1)y_{n+1}=\sum_{k=0}^n\alpha_{n-k}y_k+\sum_{k=0}^n\beta_{n-k}\lambda^ky_k+\gamma_n
\ee{yr}
for $n\ge 0$. Also assume that $\bet(0)=\beta_0\ne 0$. Then upon defining $w_n\in\IR$ by
\be{cy}
w_n=\bigg(\frac{n!}{\lambda^{n(n-1)/2}\beta_0^n}\bigg)y_n,
\ee{cy}
we have that the finite limit
$$
\lim_{n\to\infty}w_n=w_\infty
$$
exists. Further, if $w_\infty\ne 0$ then equation~\eqref{yeqq} with the initial value
$y(0)=y_0$ has no analytic solution $y(\wt t)$ in any neighborhood of $\wt t=0$.
\exx
\vv

\noindent {\bf Proof.}
Upon substituting~\eqref{cy} into~\eqref{yr}, we obtain the recursion relation
\be{cr}
w_{n+1}=\bigg(1+\frac{\alpha_0}{\lambda^n\beta_0}\bigg)w_n
+\sum_{k=0}^{n-1}\frac{\theta_n}{\theta_k}
\bigg(\beta_{n-k}+\frac{\alpha_{n-k}}{\lambda^k}\bigg)\frac{w_k}{\beta_0}+
\frac{\theta_n\gamma_n}{\beta_0},
\ee{cr}
for the terms $w_n$, where we denote
\be{the}
\theta_k=\frac{k!}{\lambda^{k(k+1)/2}\beta_0^k}.
\ee{the}
It follows from~\eqref{the2} and from~\eqref{cr} that
\be{rec1}
\begin{array}{lcl}
|w_{n+1}-w_n| &\bb \le &\bb
\ds{C\bigg|\frac{w_n}{\lambda^n\beta_0}\bigg|
+2C\sum_{k=0}^{n-1}\bigg|\frac{\theta_n\mu^n}{\theta_k\mu^k}\bigg|
\bigg|\frac{w_k}{\beta_0}\bigg|
+C\bigg|\frac{\theta_n\mu^n}{\beta_0}\bigg|}\\
\\
&\bb \le &\bb
\ds{\frac{C}{|\beta_0|}\bigg(\bigg(\frac{1}{|\lambda^n|}+2\sum_{k=0}^{n-1}\bigg|\frac{\theta_n\mu^n}{\theta_k\mu^k}\bigg|\bigg)
\max_{0\le k\le n}|w_k|+|\theta_n|\mu^n}\bigg)
\end{array}
\ee{rec1}
for every $n\ge 0$. We note that for every $k\ge 1$
\be{4y}
\bigg|\frac{\theta_k\mu^k}{\theta_{k-1}\mu^{k-1}}\bigg|=\frac{\mu k}{|\lambda^k\beta_0|}\le
\left\{\begin{array}{cl}
\ds{\frac{K_1}{2},} &\bb \quad\hbox{for }1\le k\le q,\\
\\
\ds{\frac{1}{2},} &\bb \quad\hbox{for }k>q,
\end{array}\right.
\ee{4y}
for some quantity $K_1\ge 1$ and integer $q\ge 1$.
Now fix a quantity $\del$ satisfying $|\lambda^{-1}|<\delta<1$.
Then there exists a constant $K_2\ge 1$ such that
\be{4a}
\bigg|\frac{\theta_n\mu^n}{\theta_k\mu^k}\bigg|
=\bigg|\frac{\theta_n\mu^n}{\theta_{n-1}\mu^{n-1}}\bigg|\bigg|\frac{\theta_{n-1}\mu^{n-1}}{\theta_k\mu^k}\bigg|
\le\bigg(\frac{\mu n}{|\lam^n\bet_0|}\bigg)\bigg(\frac{K_1^q}{2^{n-k-1}}\bigg)
\le\frac{K_2\del^n}{2^{n-k}},
\ee{4a}
for $0\le k\le n-1$, where $K_2$ does not depend on $n$ or $k$. With this, it follows from~\eqref{rec1} that
\be{dbnd}
\begin{array}{lcl}
|w_{n+1}-w_n| &\bb \le &\bb
\ds{\frac{C}{|\bet_0|}\bigg(\bigg(\del^n+2\sum_{k=0}^{n-1}\frac{K_2\del^n}{2^{n-k}}\bigg)\max_{0\le k\le n}|w_k|+\frac{K_2\del^n}{2^n}\bigg)}\\
\\
&\bb < &\bb
\ds{\frac{C\del^n}{|\bet_0|}\bigg(1+2K_2+K_2\bigg)\bigg(1+\max_{0\le k\le n}|w_k|\bigg)=K_3\del^n\bigg(1+\max_{0\le k\le n}|w_k|\bigg)},
\end{array}
\ee{dbnd}
where the above equality serves to define $K_3$. We claim that the sequence $w_n$ is bounded, in fact, that
\be{a1}
|w_n|\le\bigg(\prod_{k=0}^{n-1}(1+K_3\delta^k)\bigg)\bigg(|w_0|+\sum_{k=0}^{n-1}K_3\delta^k\bigg)
\ee{a1}
for every $n\ge 0$. The proof of~\eqref{a1} follows directly from~\eqref{dbnd} by induction. Indeed,~\eqref{a1}
holds for $n=0$, where the empty sum and product here are taken to be $0$ and $1$, respectively.
Assuming~\eqref{a1} holds for all for all integers from zero to $n$, we see that the right-hand side of~\eqref{a1},
with $n$ as given, is in fact an upper bound for $\ds{\max_{0\le k\le n}|w_k|}$. Therefore from~\eqref{dbnd},
$$
\begin{array}{lcl}
|w_{n+1}| &\bb < &\bb
\ds{|w_n|+K_3\delta^n\bigg(1+\max_{0\le k\le n}|w_k|\bigg)\le (1+K_3\delta^n)\max_{0\le k\le n}|w_k|+K_3\delta^n}\\
\\
&\bb \le &\bb
\ds{(1+K_3\delta^n)\bigg(\prod_{k=0}^{n-1}(1+K_3\delta^k)\bigg)\bigg(|w_0|+\sum_{k=0}^{n-1}K_3\delta^k\bigg)+K_3\delta^n}\\
\\
&\bb < &\bb
\ds{\bigg(\prod_{k=0}^{n}(1+K_3\delta^k)\bigg)\bigg(|w_0|+\sum_{k=0}^{n}K_3\delta^k\bigg),}
\end{array}
$$
as desired. This establishes~\eqref{a1} for all $n$, and thereby provides a uniform bound
$$
|w_n|<\bigg(\prod_{k=0}^{\infty}(1+K_3\delta^k)\bigg)\bigg(|w_0|+\sum_{k=0}^{\infty}K_3\delta^k\bigg)<\infty
$$
for the terms $w_n$. With this, and from~\eqref{dbnd}, it follows that $w_n$ is a Cauchy sequence
and thus the limit $w_\infty$ exists, as claimed.

Finally, to show that no analytic solution exists if $w_\infty\ne 0$, we observe by~\eqref{cy} that
$$
\lim_{n\to\infty}|y_n|^{1/n}=\lim_{n\to\infty}\bigg|\frac{\lambda^{(n-1)/2}\beta_0}{(n!)^{1/n}}\bigg||w_n|^{1/n}=\infty,
$$
where as earlier Stirling's formula is used.~\fp
\vv

\noindent {\bf Remark.}
For given coefficients $\alp(\wt t)$ and $\bet(\wt t)$ and for a given $\lam$, one sees that the quantity $w_\infty$ depends
linearly on both the initial condition $y_0$ and on the inhomogeneous term $\gam(\wt t)$. More precisely,
$$
w_\infty=w_\infty(y_0,\gam)=w_\infty^{\rm hom}y_0+w_\infty^{\rm inh}(\gam)
$$
where $w_\infty^{\rm hom}\in\IR$ is the value of $w_\infty$ that occurs with $y_0=1$ and $\gam(\wt t)\equiv 0$
identically, and where $w_\infty^{\rm inh}(\gam)$ is the value of $w_\infty$ with $y_0=0$ and $\gam$ as given.
One may regard $w_\infty^{\rm inh}:\mathcal{G}\to\IR$ as a linear map, where $\mathcal{G}$ is the vector space
(with no topology imposed) of all germs of analytic functions $\gam(\wt t)$ at $\wt t=0$.
(Of course ``hom'' and ``inh'' stand for homogeneous and inhomogeneous, respectively.)

An interesting question is whether for a given $\alp$, $\bet$, and $\lam$, with $\bet(0)\ne 0$ as in Theorem~4.2,
the map $w_\infty^{\rm inh}$ is necessarily nontrivial, that is, $w_\infty^{\rm inh}(\gam)\ne 0$ for some
$\gam\in\mathcal{G}$. We conjecture it is always nontrivial, as it would seem highly unlikely that $w_\infty^{\rm inh}(\gam)=0$
would hold for every $\gam\in\mathcal{G}$.
\vv

\noindent {\bf Remark.}
It is not clear how to extend Theorem~4.2 to the case of systems, that is,
where $y(\wt t)\in\IR^N$ in equation~\eqref{yeqq}. In such a case $\alp(\wt t)$
and $\bet(\wt t)$ are $N\times N$-matrix valued functions and $\gam(\wt t)\in\IR^N$, and
the difficulties arise from issues of commutativity. Namely, upon making the change
of variables~\eqref{cy}, where one interprets $\bet_0^n$ in the denominator as $\bet_0^{-n}$,
with $\bet_0$ assumed to be a nonsingular matrix, it is not the recursion equation~\eqref{cr} for $w_n$
that arises. Rather, one instead obtains
$$
w_{n+1}=\bigg(I+\lam^{-n}\bet_0^{-n-1}\alp_0\bet_0^n\bigg)w_n
+\sum_{k=0}^{n-1}\bet_0^{-1}\theta_n\bigg(\bet_{n-k}+\lam^{-k}\alp_{n-k}\bigg)\theta_k^{-1}w_k
+\bet_0^{-1}\theta_n\gam_n,
$$
where the various factors in the each term in the sum need not commute with one another,
and so the arguments of the above proof no longer apply.
\vv

\noindent {\bf Remark.}
In case $w_\infty=0$ in Theorem~4.2, there in fact does exist a
(necessarily unique) solution of~\eqref{yeqq} which is analytic in a neighborhood
of $t_0$; this converse to Theorem~4.2 will be shown in Theorem~4.5 below. However, there also exist other solutions
through the same initial point which are only $C^\infty$ but not analytic, and in fact there is an infinite-dimensional
space of such solutions. This is the content of Theorem~4.4.
Thus $w_\infty=0$ is no guarantee that a particular solution of interest is analytic, even though an analytic solution does exist.

An interesting open question is whether, for the original equation~\eqref{slin}
with $a(t)$, $b(t)$, $h(t)$, and $\eta(t)$ analytic and periodic, and with $w_\infty=0$
for some $y_0$ at an expansive
fixed point $t_0$ of $\eta$, there can exist a periodic solution with $y(t_0)=y_0$
which is not analytic at $t_0$.
\vv

Before proving Theorem~4.4 we require the following lemma.
By a linear space (in a Banach space), such as $E$ in this lemma, we mean simply a vector space, which need not
be closed. By an affine space, such as $E\cap\Lam^{-1}(a)$, we mean a translate of a linear space.
By an affine map we mean a map which is linear plus a constant.
\vv

\noindent {\bf Lemma~4.3.} \bxx
Let $\Lam:X\to\IR^n$ be a continuous affine map whose range is all of $\IR^n$,
where $X$ is an infinite-dimensional Banach space.
Suppose that $E\subseteq X$ is a dense linear space. Then for every $a\in\IR^n$, the set
$$
E\cap\Lam^{-1}(a)=\{x\in E\:|\:\Lam(x)=a\}
$$
is an infinite-dimensional affine space.
\exx
\vv

\noindent {\bf Proof.}
By considering the map $x\to\Lam(x)-\Lam(0)$, we may assume without loss that $\Lam$ is a linear map.
Let $\Lam_E:E\to\IR^n$ denote the restriction of $\Lam$ to $E$. Then by the continuity of $\Lam$, the
range of $\Lam_E$ is a dense subspace of $\IR^n$ and thus all of $\IR^n$. Therefore given any $a\in\IR^n$, there exists
$x_0\in E$ such that $\Lam_E(x_0)=a$. Thus it holds for any $x\in E$, that $\Lam(x)=a$ if and only if
$x-x_0\in\ker(\Lam_E)$ where $\ker(\Lam_E)=\{y\in E\:|\:\Lam_E(y)=0\}$ is the kernel of $\Lam_E$; that is,
$E\cap\Lam^{-1}(a)=x_0+\ker(\Lam_E)$. Thus it is sufficient to show that $\ker(\Lam_E)$ is infinite-dimensional.
However, the range $\IR^n$ of $\Lam_E$ is linearly isomorphic to $E/\ker(\Lam_E)$ (in the algebraic sense as vector spaces,
with no topology imposed); and as $E$ is infinite-dimensional, it follows that $\ker(\Lam_E)$ is also infinite-dimensional.~\fp
\vv

\noindent {\bf Theorem~4.4.} \bxx
Consider equation~\eqref{yeqq} with analytic coefficients $\alpha,\beta,\gamma:V\to\IR$
in some neighborhood $V\subseteq\IR$ of $\wt t=0$,
and with $|\lambda|>1$. Then there exists $\tau_0>0$ such that the following holds.
For every $y_0\in\IR$ and $\tau\in\IR$ with $0<\tau\le\tau_0$, the set
$$
\{y\in C^\infty[-\tau,\tau]\:|\:y(0)=y_0\mbox{ and }
y(\wt t)\mbox{ satisfies equation }\eqref{yeqq}\mbox{ for }|\wt t|\le|\lam|^{-1}\tau\}
$$
is an infinite-dimensional affine space. Further, all of these solutions
$y(\wt t)$ share the same (not necessarily convergent) Taylor series expansion about $\wt t=0$.
\exx
\vv

\noindent {\bf Proof.}
Assume first that $\tau_0>0$ is small enough that the coefficients $\alp$, $\bet$, and $\gam$ in~\eqref{yeqq}
are analytic for $|\wt t|\le\tau_0$; we reserve the right to reduce $\tau_0$ further. Let $K>0$ be a common bound
$$
|\alp(\wt t)|,|\bet(\wt t)|,|\gam(\wt t)|\le K,
$$
for these coefficients in this interval. Now fix $\tau$ satisfying $0<\tau\le\tau_0$.
Then any
$$
\vph\in C(I_-\cup I_+),\qquad\mbox{where}\qquad
I_-=[-\tau,-|\lam|^{-1}\tau],\qquad
I_+=[|\lam|^{-1}\tau,\tau],
$$
may be taken as an
initial condition for equation~\eqref{yeqq}, wherein the solution may be solved in
steps to the right, from $\wt t=-|\lam|^{-1}\tau$, and to the left, from $\wt t=|\lam|^{-1}\tau$, toward $\wt t=0$.
This yields a function $y(\wt t)$, continuous for $\wt t\in[-\tau,\tau]\setminus\{0\}$
and agreeing with $\vph$ on its domain $I_-\cup I_+$, and
which is $C^1$ and satisfies the differential equation~\eqref{yeqq} for $\wt t\in[-|\lam|^{-1}\tau,|\lam|^{-1}\tau]\setminus\{0\}$.
In fact this solution is $C^n$ for $\wt t\in[-|\lam|^{-n}\tau,|\lam|^{-n}\tau]\setminus\{0\}$ for every $n\ge 1$,
as one checks inductively. Let us denote this solution by $y(\wt t;\vph)$.

By a simple Gronwall argument one can show that this solution is bounded as $\wt t\to 0$, both from the left and from the right.
Indeed, upon letting
$$
\xi(\wt t)=\max_{\wt t\le |s|\le\tau}|y(s;\vph)|
$$
then one has
$$
\xi(\wt t)\le\|\vph\|+\int_{\wt t}^{|\lam|^{-1}\tau}2K\xi(s)+K\:ds
<\|\vph\|+K|\lam|^{-1}\tau+2K\int_{\wt t}^{|\lam|^{-1}\tau}\xi(s)\:ds
$$
for $0<\wt t\le|\lam|^{-1}\tau$, where $\|\cdot\|$ denotes the norm in
$C(I_-\cup I_+)$. This in turn gives
\be{ybnd}
\xi(\wt t)<(\|\vph\|+K|\lam|^{-1}\tau)e^{2K(|\lam|^{-1}\tau-\wt t)},
\ee{ybnd}
which implies that $y(\wt t;\vph)$ is bounded as $\wt t\to 0$, both from
the left and right. It follows from the differential equation
that the derivative $\dot y(\wt t;\vph)$ also remains bounded, and thus
both the left-hand and right-hand limits
$$
\Lam_-(\vph)=\lim_{\wt t\to 0-}y(\wt t;\vph)=y(0-;\vph),\qquad
\Lam_+(\vph)=\lim_{\wt t\to 0+}y(\wt t;\vph)=y(0+;\vph),
$$
exist and are finite (although they could be different). The
above formulas serve as the definitions of the two functions
$\Lam_\pm:C(I_-\cup I_+)\to\IR$, which are affine due to the linearity of the differential equation.
Also, the bound
$$
|\Lam_\pm(\vph)|\le(\|\vph\|+K|\lam|^{-1}\tau)e^{2K|\lam|^{-1}\tau}
$$
holds by~\eqref{ybnd} and implies that $\Lam_-$ and $\Lam_+$ are continuous functions.

It is further the case, as one sees by successively differentiating the differential equation, that
the derivatives have finite limits
$\ds{\lim_{\wt t\to 0\pm}}y^{(n)}(\wt t;\vph)=y^{(n)}(0\pm;\vph)$
from the left and right, for every $n\ge 1$, and thus the solution $y(\wt t;\vph)$ belongs to both
$C^n[-|\lam|^{-n}\tau,0]$ and to $C^n[0,|\lam|^{-n}\tau]$. If additionally
it is the case that $\Lam_-(\vph)=\Lam_+(\vph)$, then each pair of left- and right-hand derivatives
is equal, and so the solution belongs to $C^n[-|\lam|^{-n}\tau,|\lam|^{-n}\tau]$.

Now define $\Lam:C(I_-\cup I_+)\to\IR^2$ by
$$
\Lam(\vph)=(\Lam_-(\vph),\Lam_+(\vph)).
$$
Then the range of $\Lam$ is an affine subset of $\IR^2$, namely either a point, a straight line,
or all of $\IR^2$. We claim it is all of $\IR^2$ provided $\tau$ is small enough, and
we prove this by exhibiting a point in the interior of each of the four standard quadrants
of the plane which is in the range of $\Lam$.
To this end, let $\vph\in C(I_-\cup I_+)$ be one of the four functions with $\vph(\wt t)\equiv c_-\in\{-1,1\}$ on $I_-$
and $\vph(\wt t)\equiv c_+\in\{-1,1\}$ on $I_-$. Then $\|\vph\|=1$, and with~\eqref{ybnd} we have that
$|y(\wt t;\vph)|<(1+K|\lam|^{-1}\tau)e^{2K|\lam|^{-1}\tau}$ throughout $[-\tau,\tau]\setminus\{0\}$.
This in turn provides the bound
$$
|\dot y(\wt t;\vph)| < K\bigg(2(1+K|\lam|^{-1}\tau)e^{2K|\lam|^{-1}\tau}+1\bigg)
$$
throughout $[-|\lam|^{-1}\tau,|\lam|^{-1}\tau]\setminus\{0\}$ via the differential equation, and it follows that
$$
\begin{array}{lcl}
|\Lam_{\pm}(\vph)-c_\pm|=
|y(0\pm;\vph)-y(\pm|\lam|^{-1}\tau;\vph)|
&\bb < &\bb
K|\lam|^{-1}\tau\bigg(2(1+K|\lam|^{-1}\tau)e^{2K|\lam|^{-1}\tau}+1\bigg)\\
\\
&\bb \le &\bb
K|\lam|^{-1}\tau_0\bigg(2(1+K|\lam|^{-1}\tau_0)e^{2K|\lam|^{-1}\tau_0}+1\bigg).
\end{array}
$$
Therefore, if $\tau_0$ is chosen small enough that
$$
K|\lam|^{-1}\tau_0\bigg(2(1+K|\lam|^{-1}\tau_0)e^{2K|\lam|^{-1}\tau_0}+1\bigg)\le 1
$$
holds, then $\Lam_{\pm}(\vph)\ne 0$ and this quantity has the same sign as $c_\pm$. In particular,
for these four choices of $\vph$, the points $\Lam(\vph)\in\IR^2$ lie in the interiors of the four
standard quadrants of the plane, as desired, and the range of $\Lam$ is all of $\IR^2$.

Now let
$$
\begin{array}{lcr}
E &\bb = &\bb
\{\vph\in C^\infty(I_-\cup I_+)\:|\:
\psi_1^{(n)}(\pm|\lam|^{-1}\tau)=\psi_2^{(n)}(\pm|\lam|^{-1}\tau)\mbox{ for every }n\ge 0,\\
\\
&\bb &\bb
\mbox{ where }\psi_1(\wt t)=\dot\vph(\wt t)-\alp(\wt t)\vph(\wt t)-\gam(\wt t)\mbox{ and }
\psi_2(\wt t)=\bet(\wt t)\vph(\lam \wt t)\}.
\end{array}
$$
One checks that the above solution $y(\wt t;\vph)$ through any $\vph\in E$ is $C^\infty$ in
neighborhoods of $\wt t=\pm|\lam|^{-1}\tau$ due to the compatibility conditions which define $E$. Thus,
by induction, $y(\wt t;\vph)$ is $C^\infty$ throughout the set $[-\tau,\tau]\setminus\{0\}$, and as noted,
it is in fact $C^\infty$ throughout $[-\tau,\tau]$ if $\Lam_-(\vph)=\Lam_+(\vph)$.

The set $E$ is a linear space which is dense in $C(I_-\cup I_+)$,
and so by Lemma~4.3, for any $a=(a_-,a_+)\in\IR^2$ the set
$\{\vph\in E\:|\:\Lam(\vph)=a\}$ is an infinite-dimensional affine space.
In particular, if we take $a_-=a_+=y_0$ for any given $y_0\in\IR$, then
the elements of $\{\vph\in E\:|\:\Lam(\vph)=(y_0,y_0)\}$ give the
desired $C^\infty$ solutions $y(\wt t;\vph)$ of equation~\eqref{yeqq}
on the interval $[-\tau,\tau]$.

The final sentence in the statement of the theorem, concerning the Taylor series,
follows directly by repeatedly differentiating equation~\eqref{yeqq}.~\fp
\vv

The next result provides a converse of sorts to Theorem~4.2.
\vv

\noindent {\bf Theorem~4.5.} \bxx
Consider the setting of Theorem~4.2, and suppose $y_0$ is such that $w_\infty=0$. Then
there exist constants $A>0$ and $\nu>0$ such that
\be{rd2}
|y_n|\le A\nu^n
\ee{rd2}
for every $n\ge 0$. The quantities $y_n$ are therefore the coefficients of
a solution of equation~\eqref{yeqq}, with $y(0)=y_0$,
which is analytic for $\wt t$ in a neighborhood of $\wt t=0$.
\exx
\vv

\noindent {\bf Proof.}
We shall show that there exist $A>0$ and $\nu>0$ such that for every integer $n\ge 0$, we have the bound
\be{nminq}
|w_n|\le A\nu^m\bigg|\frac{\lam^n\theta_n}{\lam^{n-m}\theta_{n-m}}\bigg|\qquad\mbox{whenever }n\ge m\ge 0,
\ee{nminq}
on the quantities given by the recursion~\eqref{cr}. Thus taking $m=n$ gives the bound
$|w_n|\le A\nu^n|\lam^n\theta_n|$ which is equivalent to the desired bound~\eqref{rd2}.
As before, we assume the bounds~\eqref{the2} on the coefficients for some $C>0$ and $\mu>0$, and the notation~\eqref{the}.

The inequality~\eqref{nminq} will be established by an induction on $m$.
Assuming that $\ds{\lim_{n\to\infty}}w_n=0$,
let us first observe, from~\eqref{the2} and~\eqref{cr}, that for every $n\ge 0$
\be{wn}
\begin{array}{lcl}
\ds{|w_n|\le\sum_{k=n}^\infty |w_{k+1}-w_k|}
&\bb \le &\bb
\ds{\sum_{k=n}^\infty
\bigg|\frac{\alp_0w_k}{\lam^k\bet_0}+\sum_{j=0}^{k-1}\frac{\theta_k}{\theta_j}
\bigg(\beta_{k-j}+\frac{\alpha_{k-j}}{\lambda^j}\bigg)\frac{w_j}{\beta_0}+
\frac{\theta_k\gamma_k}{\beta_0}\bigg|}\\
\\
&\bb \le &\bb
\ds{\frac{C}{|\bet_0|}\sum_{k=n}^\infty\bigg(\bigg|\frac{w_k}{\lam^k}\bigg|+2\sum_{j=0}^{k-1}
\bigg|\frac{\theta_k\mu^k}{\theta_j\mu^j}\bigg||w_j|+|\theta_k\mu^k|\bigg)}.
\end{array}
\ee{wn}
As $w_n$ is a bounded sequence, we
may fix $A\ge 1$ such that
$|w_n|\le A$
for every $n\ge 0$. Let such $A$ be fixed for the remainder of the proof.
This establishes~\eqref{nminq} for $m=0$ and so begins the induction.
The quantity $\nu$ will be chosen later, but of course it must be independent of $m$ and $n$.

We now proceed with the inductive step from $m-1$ to $m$.
Assume for some $m\ge 1$ that
\be{indh}
|w_n|\le A\nu^i\bigg|\frac{\lam^n\theta_n}{\lam^{n-i}\theta_{n-i}}\bigg|\qquad\mbox{whenever }n\ge i\ge 0\mbox{ and }0\le i\le m-1.
\ee{indh}
That is, assume the inductive hypothesis up to $m-1$.
We must establish~\eqref{nminq} for $m$, in the range $n\ge m$, as shown.
For such $n$ we have from~\eqref{wn} that
\be{wn2}
\begin{array}{lcl}
|w_n| &\bb \le &\bb
\ds{\frac{C}{|\bet_0|}\sum_{k=n}^\infty\bigg(
\bigg|\frac{w_k}{\lam^k}\bigg|+
2\sum_{j=0}^{m-2}\bigg|\frac{\theta_k}{\theta_j}\bigg|\mu^{k-j}|w_j|
+2\sum_{j=m-1}^{k-1}\bigg|\frac{\theta_k}{\theta_j}\bigg|\mu^{k-j}|w_j|
+|\theta_k|\mu^k\bigg)}\\
\\
&\bb \le &\bb
\ds{\frac{AC}{|\bet_0|}\sum_{k=n}^\infty\bigg(
\nu^{m-1}\bigg|\frac{\theta_k}{\lam^{k-m+1}\theta_{k-m+1}}\bigg|+
2\sum_{j=0}^{m-2}\bigg|\frac{\theta_k}{\theta_j}\bigg|\mu^{k-j}\nu^j|\lam^j\theta_j|}\\
\\
&\bb &\bb
\qquad\qquad
\ds{+2\sum_{j=m-1}^{k-1}\bigg|\frac{\theta_k}{\theta_j}\bigg|\mu^{k-j}\nu^{m-1}\bigg|\frac{\lam^j\theta_j}{\lam^{j-m+1}\theta_{j-m+1}}\bigg|
+|\theta_k|\mu^k\bigg)}\\
\\
&\bb = &\bb
\ds{\frac{AC}{|\bet_0|}\sum_{k=n}^\infty\bigg(
\bigg|\frac{\theta_k}{\lam^{k-m+1}\theta_{k-m+1}}\bigg|\nu^{m-1}+
2\sum_{j=0}^{m-2}|\lam^j\theta_k|\mu^{k-j}\nu^j}\\
\\
&\bb &\bb
\qquad\qquad\ds{+2\sum_{j=m-1}^{k-1}\bigg|\frac{\theta_k}{\lam^{-m+1}\theta_{j-m+1}}\bigg|\mu^{k-j}\nu^{m-1}
+|\theta_k|\mu^k\bigg)}.
\end{array}
\ee{wn2}
Note in particular that for the range $0\le j\le m-2$, in estimating $|w_j|$
we have used the induction hypothesis~\eqref{indh} with $i=j$, while for $m-1\le j\le k$ we have taken $i=m-1$.
Note that the assumption $A\ge 1$ was used in dealing with the final term $|\theta_k|\mu^k$.
Also, when $m=1$ then the sum from $j=0$ to $j=m-2$ by convention is zero.

We shall impose the requirement on the quantity $\nu$ (yet to be chosen) that
$\nu\ge 2|\lam|^{-1}\mu$,
and with this we have the estimate
$$
\begin{array}{lcl}
\ds{2\sum_{j=0}^{m-2}|\lam^j\theta_k|\mu^{k-j}\nu^j}
&\bb < &\bb
\ds{2|\lam^{m-1}\theta_k|\bigg(\frac{\mu^{k-m+1}\nu^{m-1}}{|\lam|\mu^{-1}\nu-1}\bigg)}\\
\\
&\bb \le &\bb
\ds{2|\lam^{m-1}\theta_k|\mu^{k-m+1}\nu^{m-1}
=2\bigg|\frac{\theta_k}{\lam^{-m+1}\theta_0}\bigg|\mu^{k-m+1}\nu^{m-1}}.
\end{array}
$$
This quantity may now be incorporated into the final sum in~\eqref{wn2} as the $j=m-1$ term.
Additionally, we have that
$$
|\theta_k|\mu^k=\bigg|\frac{\theta_k}{\theta_0}\bigg|\mu^k
\le \bigg|\frac{\theta_k}{\lam^{-m+1}\theta_0}\bigg|\mu^{k-m+1}\nu^{m-1},
$$
so this term also may be incorporated into the $j=m-1$ term. With this we obtain
$$
\begin{array}{lcl}
|w_n| &\bb \le &\bb
\ds{\frac{AC}{|\bet_0|}\sum_{k=n}^\infty\bigg(
\bigg|\frac{\theta_k}{\lam^{k-m+1}\theta_{k-m+1}}\bigg|\nu^{m-1}
+2\sum_{j=0}^{m-2}|\lam^j\theta_k|\mu^{k-j}\nu^j}\\
\\
&\bb &\bb
\qquad\qquad\ds{+2\sum_{j=m-1}^{k-1}\bigg|\frac{\theta_k}{\lam^{-m+1}\theta_{j-m+1}}\bigg|\mu^{k-j}\nu^{m-1}
+|\theta_k|\mu^k\bigg)}\\
\\
&\bb < &\bb
\ds{\frac{AC}{|\bet_0|}\sum_{k=n}^\infty\bigg(
\bigg|\frac{\theta_k}{\lam^{k-m+1}\theta_{k-m+1}}\bigg|\nu^{m-1}
+5\sum_{j=m-1}^{k-1}\bigg|\frac{\theta_k}{\lam^{-m+1}\theta_{j-m+1}}\bigg|\mu^{k-j}\nu^{m-1}\bigg)}.
\end{array}
$$
We require that this quantity be bounded above by the right-hand side of~\eqref{nminq}, that is, we require that
\be{5th}
\frac{C}{|\lam\bet_0|}\sum_{k=n}^\infty\bigg(
\bigg|\frac{\theta_k\theta_{n-m}}{\lam^k\theta_{k-m+1}\theta_n}\bigg|
+5\sum_{j=m-1}^{k-1}
\bigg|\frac{\theta_k\theta_{n-m}}{\theta_{j-m+1}\theta_n}\bigg|\mu^{k-j}\bigg)
\le\nu
\ee{5th}
hold for every $n\ge m$. As such we seek an upper bound for $|\theta_k\theta_{n-m}|/|\theta_{j-m+1}\theta_n|$,
and to this end we first consider $|\theta_k\theta_{n-m}|/|\theta_{k-m}\theta_n|$,
where $k\ge n\ge m$. From~\eqref{the} we obtain
$$
\frac{\theta_k\theta_{n-m}}{\theta_{k-m}\theta_n}=\bigg(\frac{k!(n-m)!}{(k-m)!n!}\bigg)\bigg(\frac{1}{\lam^{(k-n)m}}\bigg)
$$
after a short calculation. Next observe that
\be{x3}
\frac{k!(n-m)!}{(k-m)!n!}
=\bigg(\frac{k!}{(k-m)!}\bigg)\bigg(\frac{(n-m)!}{n!}\bigg)=\prod_{i=0}^{m-1}\bigg(\frac{k-i}{n-i}\bigg)\le(k-n+1)^m,
\ee{x3}
and also that
\be{x2}
\frac{k!(n-m)!}{(k-m)!n!}
=\bigg(\frac{k!}{n!}\bigg)\bigg(\frac{(n-m)!}{(k-m)!}\bigg)
=\prod_{i=n+1}^{k}\bigg(\frac{i}{i-m}\bigg)\le(m+1)^{k-n}.
\ee{x2}
Both of these bounds will be used,
the choice depending on the range of the indices.
We conclude from~\eqref{x3} and~\eqref{x2} that
$$
\bigg|\frac{\theta_k\theta_{n-m}}{\theta_{j-m+1}\theta_n}\bigg|\mu^{k-j}
\le
\min\{(k-n+1)^m,(m+1)^{k-n}\}\bigg(\frac{\mu}{|\lam|^{(k-n)m}}\bigg)
\bigg|\frac{\theta_{k-m}\mu^{k-m}}{\theta_{j-m+1}\mu^{j-m+1}}\bigg|,
$$
and therefore
\be{4b}
\begin{array}{l}
\ds{\bigg|\frac{\theta_k\theta_{n-m}}{\lam^k\theta_{k-m+1}\theta_n}\bigg|
+5\sum_{j=m-1}^{k-1}\bigg|\frac{\theta_k\theta_{n-m}}{\theta_{j-m+1}\theta_n}\bigg|\mu^{k-j}}\\
\\
\qquad\qquad\ds{\le\min\{(k-n+1)^m,(m+1)^{k-n}\}}\\
\\
\qquad\qquad\quad\ds{\times\bigg(\frac{\mu}{|\lam|^{(k-n)m}}\bigg)
\bigg(\bigg|\frac{\theta_{k-m}}{\lam^k\theta_{k-m+1}\mu}\bigg|
+5\sum_{j=0}^{k-m}\bigg|\frac{\theta_{k-m}\mu^{k-m}}{\theta_{j}\mu^{j}}\bigg|\bigg).}
\end{array}
\ee{4b}
(Note that we have shifted the index $j$ in the final sum.)
Now with $K_2$ and $\del$ as in~\eqref{4y} and~\eqref{4a}
in the proof of Theorem~4.2, with $|\lam^{-1}|<\del<1$, we have from~\eqref{4a} and~\eqref{4b} that
$$
\begin{array}{l}
\ds{\bigg|\frac{\theta_k\theta_{n-m}}{\lam^k\theta_{k-m+1}\theta_n}\bigg|
+5\sum_{j=m-1}^{k-1}\bigg|\frac{\theta_k\theta_{n-m}}{\theta_{j-m+1}\theta_n}\bigg|\mu^{k-j}}\\
\\
\qquad\qquad
\ds{\le\min\{(k-n+1)^m,(m+1)^{k-n}\}}\\
\\
\qquad\qquad\quad
\ds{\times\bigg(\frac{\mu}{|\lam|^{(k-n)m}}\bigg)\bigg(\bigg|\frac{\bet_0}{\lam^{m-1}(k-m+1)\mu}\bigg|
+5\sum_{j=0}^{k-m-1}\frac{K_2\del^{k-m}}{2^{k-m-j}}+5\bigg)}\\
\\
\qquad\qquad
\ds{<\min\{(k-n+1)^m,(m+1)^{k-n}\}\bigg(\frac{\mu}{|\lam|^{(k-n)m}}\bigg)\bigg(\bigg|\frac{\bet_0}{\mu}\bigg|+5K_2\del^{k-m}+5\bigg)}\\
\\
\qquad\qquad
\ds{\le K_4\min\{(k-n+1)^m,(m+1)^{k-n}\}\bigg(\frac{1}{|\lam|^{(k-n)m}}\bigg)},
\end{array}
$$
where $K_4=|\bet_0|+5\mu(K_2+1)$.
Also, there exist a quantity $K_5>0$ and an integer $c\ge 0$ such that
$$
\frac{k+1}{|\lam|^k}\le K_5\mbox{ for every }k\ge 0,\qquad
\frac{k+1}{|\lam|^k}\le \del^k\mbox{ for every }k\ge c.
$$
Therefore
$$
\begin{array}{l}
\ds{\sum_{k=n}^\infty
\bigg(\bigg|\frac{\theta_k\theta_{n-m}}{\lam^k\theta_{k-m+1}\theta_n}\bigg|
+5\sum_{j=m-1}^{k-1}\bigg|\frac{\theta_k\theta_{n-m}}{\theta_{j-m+1}\theta_n}\bigg|\mu^{k-j}\bigg)}\\
\\
\qquad\qquad\ds{<K_4\sum_{k=0}^\infty\min\{(k+1)^m,(m+1)^k\}\bigg(\frac{1}{|\lam|^{km}}\bigg)}\\
\\
\qquad\qquad\ds{\le K_4\sum_{k=0}^{c-1}\bigg(\frac{m+1}{|\lam|^{m}}\bigg)^k+K_4\sum_{k=c}^\infty\bigg(\frac{k+1}{|\lam|^{k}}\bigg)^m}\\
\\
\qquad\qquad\ds{\le K_4\sum_{k=0}^{c-1}K_5^k+K_4\sum_{k=c}^\infty\del^{km}
\le K_4\sum_{k=0}^{c-1}K_5^k+\frac{K_4\del^c}{1-\del}=K_6},
\end{array}
$$
where the final equality serves as the definition of the quantity $K_6$. It follows now
that the required condition~\eqref{5th} holds as long as $\nu$ is large enough that $\nu\ge CK_6|\lam\bet_0|^{-1}$.
With the other required condition $\nu\ge 2|\lam|^{-1}\mu$ also holding, the proof of the theorem is complete.~\fp

\section{Coexistence of Analyticity and Nonanalyticity}

One expects that typically (generically) $w_\infty\ne 0$, and that the
case of $w_\infty=0$ is exceptional, of codimension one in some sense.
Here we shall present a class of examples which are linear equations with periodic coefficients,
for which $w_\infty\ne 0$ is realized for a specific periodic solution. In fact, we shall
show there are periodic solutions for which $\eta$ possesses both a contractive fixed
point and an expansive one with $w_\infty\ne 0$, which we may
term {\bf coexistence} of analyticity and nonanalyticity. As such, the
solution is analytic for some values of $t$ and not analytic for others.
Equivalently, both the sets $\mathcal{A}$ and $\mathcal{N}$ in~\eqref{an}
are nonempty for such a solution.

We consider a class of integral equations of the form
\be{inteq}
\kappa x(t)=\int_{t-r(t)}^t \rho(s)x(s)\:ds,
\ee{inteq}
where we assume both analyticity and a periodicity condition, namely, that
\be{v1}
\begin{array}{l}
r,\rho:\IR\to\IR\mbox{ are analytic, and}\\
\\
r(t+2\pi)=r(t),\qquad\rho(t+2\pi)=\rho(t),\qquad\mbox{for every }t\in\IR.
\end{array}
\ee{v1}
Here $\kappa$ is a parameter (an eigenvalue) to be found along with the solution $x(t)$.
Generally, we shall be interested in $2\pi$-periodic solutions, that is, solutions
for which $x(t+2\pi)=x(t)$ for every $t\in\IR$. For any such solution for which $\kappa\ne 0$, we obtain
the differential equation
$$
\kappa\dot x(t)=\rho(t)x(t)-(1-\dot r(t))\rho(t-r(t))x(t-r(t))
$$
upon differentiating~\eqref{inteq}.
Note that any such periodic solution is $C^\infty$ everywhere.

If $x(t)$ is a periodic solution as above, with~\eqref{v1} holding, we observe
that it also satisfies the modified equation
\be{u++}
\kappa\dot x(t)=\rho(t)x(t)-\dot\eta(t)\rho(\eta(t))x(\eta(t)),\qquad
\eta(t)=t-r(t)+2\pi m,
\ee{u++}
for any integer $m$. In fact it is to this modified equation that we shall apply our theorems.

We begin by providing, in Theorem~5.1 below, a general existence-uniqueness result for
a broad class of periodic integral equations~\eqref{inteq}. The proof of this
theorem, which does not require analyticity, involves the Kre\u\i n--Rutman Theorem (see~\cite{kr})
and standard arguments, and we provide a proof for the reader's convenience.
We also refer the reader to~\cite{nkr}, which gives a detailed analysis of linear operators
very closely related to equation~\eqref{inteq} using similar arguments.

We note in passing that even under the given hypotheses of Theorem~5.1, there may exist
other $2\pi$-periodic solutions of equation~\eqref{inteq} which take both positive and negative values.
\vv

\noindent {\bf Theorem~5.1.} \bxx
Consider equation~\eqref{inteq} where $r:\IR\to(0,\infty)$ and $\rho:\IR\to(0,\infty)$
are continuous and positive functions, and satisfy $r(t+2\pi)=r(t)$ and $\rho(t+2\pi)=\rho(t)$
for every $t\in\IR$. Then there exists a unique $\kappa>0$ such that equation~\eqref{inteq}
possesses a nontrivial nonnegative $2\pi$-periodic solution, and moreover, this solution is strictly
positive; that is,
$x(t+2\pi)=x(t)>0$ for every $t\in\IR$.
Further, this solution is unique up to scalar multiple and we have the bounds
\be{radbnd}
\min_{t\in\IR}\bigg(\int^t_{t-r(t)}\rho(s)\:ds\bigg)\le\kappa\le\max_{t\in\IR}\bigg(\int^t_{t-r(t)}\rho(s)\:ds\bigg)
\ee{radbnd}
for the eigenvalue $\kappa$.
\exx
\vv

The following theorem is a main result of this section, providing a class of
examples for which $w_\infty\ne 0$ holds.
\vv

\noindent {\bf Theorem~5.2.} \bxx
Consider the equation
\be{ie}
\kappa x(t)=\int^t_{t-r(t)}x(s)\:ds,\qquad
r(t)=-(\lam-1)\sin t+2\pi m,
\ee{ie}
where $1<\lam<2\pi m+1$, with $m$ an integer. Then there exists $\lam_*$, independent of $m$, such that the following holds.
If $\lam\ge \lam_*$ and $m$ also satisfies $2\pi m\le 2\lam+1$,
and $x(t)$ is the unique positive $2\pi$-periodic solution given by Theorem~5.1 with $\rho(t)\equiv 1$ identically,
then $x(t)$ is not analytic in any neighborhood of $t=0$. In particular, $w_\infty\ne 0$ for the quantity
in Theorem~4.2, for the equation~\eqref{yeqq} obtained from~\eqref{u++}
by the transformation~\eqref{con-}, \eqref{cov} with $t_0=0$.
\exx
\vv

Before proving the above results, let us observe that in the analytic case
of Theorem~5.1, and where $\eta$ possesses an expansive fixed point $t_0$, it is possible
for $x(t)$ to be analytic for every $t$ (and thus $w_\infty=0$ for the point $t_0$). Indeed,
let $r:\IR\to(0,\infty)$ and $z:\IR\to(0,\infty)$ both be $2\pi$-periodic and analytic,
and suppose that
$$
r(t_0)=2\pi m,\qquad
|1-\dot r(t_0)|>1,
$$
for some $t_0$ and positive integer $m$. Set
$$
x(t)=\int^t_{t-r(t)}z(s)\:ds,\qquad
\rho(t)=\frac{z(t)}{x(t)}.
$$
Then $x(t)$ is the solution given in Theorem~5.1, with $\kappa=1$,
and it thus satisfies equation~\eqref{u++};
and certainly $x(t)$ is analytic, and $t_0$ is an expansive fixed point of this equation with $w_\infty=0$.

For the reader's convenience, in preparation for the proof of
Theorem~5.1 we give a statement of the Kre\u\i n--Rutman Theorem,
and we remark upon certain subtleties which are sometimes overlooked.
Let $X$ be a real Banach space. By a {\bf cone} in $X$ we mean a convex
set $X^+\subseteq X$ such that $\tau X^+\subseteq X^+$ for every $\tau\ge 0$,
and $X^+\cap (-X^+)=\{0\}$. (Generally we shall consider closed cones, namely
cones which are closed sets.) A cone is called {\bf total} if $S$ is dense in $X$
where $S=\{x-y\:|\:x,y\in X^+\}$, and it is called {\bf reproducing} if $S=X$.
In infinite dimensions it may easily happen that a closed cone is total but not reproducing.
The Kre\u\i n--Rutman Theorem states that if $L:X\to X$ is a compact linear operator such that
$L(X^+)\subseteq X^+$ where $X^+\subseteq X$ is a closed total cone, and if also $\rad(L)>0$
where $\rad(L)$ denotes the spectral radius of $L$, then there exists $x\in X^+\setminus\{0\}$
such that $Lx=\rad(L)x$. The assumption that $\rad(L)>0$ is crucial in infinite dimensions.
Also, without further assumptions, the eigenvector $x$ need not be unique (up to scalar multiple).

The Kre\u\i n--Rutman Theorem has been generalized in a variety of directions by many authors.
We refer to~\cite{genkr} for some references to the extensive literature on the subject.
\vv

\noindent {\bf Proof of Theorem~5.1.}
Define an operator $L:X\to X$ by
$$
(Lx)(t)=\int_{t-r(t)}^t \rho(s)x(s)\:ds
$$
on the Banach space
$$
X=\{x:\IR\to\IR\:|\:x(t)\mbox{ is continuous, and }x(t+2\pi)=x(t)\mbox{ for all }t\in\IR\}
$$
endowed with the supremum norm. Then~\eqref{inteq} may be written as
\be{opeq}
Lx=\kappa x.
\ee{opeq}
The operator $L$ is compact (due to the Ascoli--Arzel\`a Theorem)
and is positive with respect to the closed reproducing cone
$$
X^+=\{x\in X\:|\:x(t)\ge 0\mbox{ for every }t\in\IR\}
$$
of nonnegative functions in $X$.
Thus by the Kre\u\i n--Rutman Theorem, there exists $x_0\in X^+\setminus\{0\}$
for which~\eqref{opeq} holds with $x=x_0$ and $\kappa=\rad(L)$, provided that $\rad(L)>0$.

To prove that $\rad(L)>0$ we shall establish the inequalities~\eqref{radbnd} with $\rad(L)$ in place of $\kappa$.
For any $x\in X$ we have that
$$
\|Lx\|=\max_{t\in\IR}\bigg|\int_{t-r(t)}^t \rho(s)x(s)\:ds\bigg|
\le \max_{t\in\IR}\bigg(\int_{t-r(t)}^t \rho(s)\:ds\bigg)\|x\|
$$
and so
$$
\rad(L)\le\|L\|\le\max_{t\in\IR}\bigg(\int_{t-r(t)}^t \rho(s)\:ds\bigg).
$$
On the other hand, letting $e\in X^+$ denote the constant function with $e(t)\equiv 1$ identically, we have that
\be{einq}
K^ne\le L^ne,\qquad\mbox{where}\qquad
K=\min_{t\in\IR}\bigg(\int_{t-r(t)}^t \rho(s)\:ds\bigg),
\ee{einq}
for every $n\ge 1$, where $\le$ denotes the partial ordering in $X$ given by the cone $X^+$.
Indeed, the inequality~\eqref{einq} is clear for $n=1$, and is easily proved for $n>1$ by inducting on $n$.
Taking norms and then $n^{\rm th}$ roots now gives
$K=\|K^ne\|^{1/n}\le\|L^ne\|^{1/n}\le\|L^n\|^{1/n}$, and letting $n\to\infty$ gives
$$
\min_{t\in\IR}\bigg(\int_{t-r(t)}^t \rho(s)\:ds\bigg)\le\rad(L).
$$
This establishes the inequalities~\eqref{radbnd} with $\rad(L)$ in place
of $\kappa$, and also establishes the existence of an eigenvector $x_0\in X^+\setminus\{0\}$ with eigenvalue $\rad(L)>0$.

Further arguments are needed to show that $x_0$ is the unique eigenvector in $X^+\setminus\{0\}$,
and the strict positivity of $r(t)$ and $\rho(t)$ is crucial here.
Note first that there exists $m\ge 1$ such that $L^m$ maps $X^+\setminus\{0\}$ into the interior
$$
\interior(X^+)=\{x\in X\:|\:x(t)>0\mbox{ for every }t\in\IR\}
$$
of the cone. Indeed, if $\ds{r_0=\min_{t\in\IR}r(t)}$, then whenever $x(t_0)>0$ for some
$x\in X^+$ and $t_0\in\IR$, we have $(Lx)(t)>0$
for every $t\in[t_0,t_0+r_0]$. Now with $r_0>0$, fix $m$ so that $mr_0\ge 2\pi$; then $L^mx\in\interior(X^+)$ for every $x\in X^+\setminus\{0\}$.
This immediately implies that for any eigenvector in the cone, say for $Lx=\kappa x$ with $x\in X^+\setminus\{0\}$ and $\kappa\ge 0$,
we have that $x\in\interior(X^+)$ and $\kappa>0$.

Suppose now that in addition to the eigenvector $x_0$ with eigenvalue $\rad(L)$
obtained above, we have another eigenvector
$y\in\interior(X^+)$ with eigenvalue $\kappa>0$. Define quantities
\be{laminq}
\tau_1=\sup\{\tau\ge 0\:|\:\tau x_0\le y\},\qquad
\tau_2=\sup\{\tau\ge 0\:|\:\tau y\le x_0\},
\ee{laminq}
both of which are positive. Then $\tau_1x_0\le y$ and so by applying $L$ we have that
$\tau_1\rad(L)x_0\le\kappa y$. Thus $\tau_1\rad(L)\kappa^{-1}\le\tau_1$ by the definition of $\tau_1$, and so $\rad(L)\le\kappa$.
A similar argument with the second equation in~\eqref{laminq} yields the opposite inequality, and thus
$\kappa=\rad(L)$. Now consider $z=y-\tau_1x_0$. It is enough to show that $z=0$, the zero element of $X$, in order to establish
the uniqueness of $x_0$. Certainly $Lz=\rad(L)z$ and $z\in X^+$; thus by the remarks above,
if $z\ne 0$ then $z\in\interior(X^+)$. However, it is clear
from the definition of $\tau_1$ that $z(t)=0$ for some $t\in\IR$ and thus $z\not\in\interior(X^+)$.
We conclude that $z=0$, and with this the theorem is proved.~\fp
\vv

In proving Theorem~5.2, we use the conjugacy $\sigma(t)$ given by Lemma~4.1 applied to
the modified equation~\eqref{u++}, where as noted
\be{eta}
\eta(t)=t-r(t)+2\pi m=t+(\lam-1)\sin t.
\ee{eta}
Clearly this conjugacy depends on the choice of $\lam$, and so it is necessary to show
that it is well-behaved for large $\lam$.
Interestingly, $\sig$ becomes better-behaved as $\lam\to\infty$ in the sense that
$\sig(t)\to t$ uniformly on the disc $|t|\le 1$. The following lemma addresses this point.
Its proof follows the broad outlines of the proof of Lemma~4.1 in~\cite{cargam}, but it also provides
some explicit estimates which will play an important role in the proof of Theorem~5.2.
\vv

\noindent {\bf Lemma~5.3.} \bxx
Let $\sig(t)$ be as in Lemma~4.1, where $\eta(t)$ is as in~\eqref{eta} with $\lam>1$ and $t_0=0$.
Then there exists $\lam_{**}>0$ such that
\be{qbn}
|\sigma(t)-t|\le \frac{|t|^2}{\lam},\qquad\hbox{for }|t|\le 1,
\ee{qbn}
provided $\lam\ge \lam_{**}$. Also, we have the bound
\be{sigb}
|\sigma_n|\le \frac{1}{\lam},\qquad\hbox{for }n\ge 2,
\ee{sigb}
for the coefficients of the Taylor series
\be{sg}
\sigma(t)=\sum_{n=1}^\infty\sigma_nt^n,
\ee{sg}
for such $\lam$.
\exx
\vv

\noindent {\bf Proof.}
As in Lemma~4.1 we have that
\be{sl}
\sigma(t)+(\lambda-1)\sin\sigma(t)=\sigma(\lambda t)
\ee{sl}
for $t$ in a neighborhood of zero, following the final equation in~\eqref{con-}.
Let us denote $g(u)=u^{-1}\sin u-1$, which is an entire function of $u$ with $g(0)=0$, and as well denote
$\del=\lambda^{-1}$. Also let us write
$\sigma(t)=t+t\zeta(t)$ where we require that $\zeta(t)$ be analytic in a neighborhood of
zero with $\zeta(0)=0$. Then after some manipulations, including replacing $t$ by $\del t$,
equation~\eqref{sl} takes the form
\be{t1}
\zeta(t)=\zeta(\del t)+(1-\del)(1+\zeta(\del t))g(\del t+\del t\zeta(\del t)).
\ee{t1}
We shall in fact obtain a solution to this equation via a fixed-point theorem,
for sufficiently small $\del$. Define the Banach space
$$
Z=\{\zeta:\overline{D_1(0)}\to\IC\:|\:\zeta(\cdot)\hbox{ is analytic in }D_1(0)\hbox{ and continuous in }\overline{D_1(0)},
\hbox{ with }\zeta(0)=0\},
$$
where as earlier $D_1(0)=\{t\in\IC\:|\:|t|<1\}$, and where we take the supremum norm for elements of $Z$.
Then the right-hand side of~\eqref{t1} defines a continuous nonlinear operator $\mathcal{T}:Z\to Z$ by
$$
(\mathcal{T}\zeta)(t)=\zeta(\del t)+(1-\del)(1+\zeta(\del t))g(\del t+\del t\zeta(\del t))
$$
for $|t|\le 1$, as long as $0<\del\le 1$.
(In fact, $(\mathcal{T}\zet)(t)$ is analytic for $|t|<\del^{-1}$ and continuous for $|t|\le\del^{-1}$.)
We claim that if additionally $\del<1$ then $\mathcal{T}$ is compact, as the following
argument shows.
Let such $\del$ be fixed and let $B_\eps=\{\zeta\in Z\:|\:\|\zeta\|\le \eps\}$ denote the closed ball of radius $\eps$ in $Z$. Also let
$\ds{G(\eps)=\sup_{|u|\le \eps}|g(u)|}$, noting that $G(\eps)$ is a continuous function with $G(0)=0$. Then
by the Schwarz maximum principle
$$
|\zeta(t)|\le\|\zeta\||t|\le \eps|t|,\qquad\hbox{for }|t|\le 1,
$$
for any $\zeta\in B_\eps$, and so for such $\zeta$
$$
\begin{array}{lcl}
|(\mathcal{T}\zeta)(t)| &\bb \le &\bb
\eps\del |t|+(1-\del)(1+\eps\del |t|)G(\del |t|(1+\eps\del |t|))\\
\\
&\bb \le &\bb
\left\{\begin{array}{ll}
\eps+(1-\del)(1+\eps)G(1+\eps)=\eps_0, &\bb \quad \hbox{for }|t|\le\del^{-1},\\
\\
\eps\del +(1-\del)(1+\eps\del )G(\del+\eps\del^2)=\eps_1, &\bb \quad \hbox{for }|t|\le 1,
\end{array}\right.
\end{array}
$$
where the above formulas serve as the definitions of $\eps_i=\eps_i(\eps)$ for $i=0,1$.
Note in particular that $\mathcal{T}\zeta\in B_{\eps_1}$. A bound for
the derivative of $(\mathcal{T}\zeta)(t)$, for $|t|\le 1$, is obtained from the Cauchy integral
formula, namely
$$
\bigg|\frac{d}{dt}(\mathcal{T}\zeta)(t)\bigg|=\bigg|\frac{1}{2\pi i}\int_{s=\del^{-1}}\frac{(\mathcal{T}\zeta)(s)}{(s-t)^2}\:ds\bigg|
\le\frac{\eps_0\del^{-1}}{(\del^{-1}-1)^2},
$$
and so the image of $B_\eps$ under $\mathcal{T}$ is an equicontinuous set.
Thus $\mathcal{T}$ is a compact map from $B_\eps$ into $B_{\eps_1}$. If in fact $\eps_1\le \eps$, that is, if
$\eps\del +(1-\del)(1+\eps\del)G(\del+\eps\del^2)\le \eps$, or equivalently, if
\be{rcond}
(1+\eps\del)G(\del+\eps\del^2)\le \eps,
\ee{rcond}
then by the Schauder Fixed Point Theorem there exists a fixed point of $\mathcal{T}$ in $B_\eps$.
Such a fixed point satisfies the desired equation~\eqref{t1} and provides the bound
$|\sigma(t)-t|=|t\zeta(t)|\le \eps|t|^2$ for the conjugacy. The desired bound~\eqref{qbn} is obtained
by taking $\eps=\del=\lambda^{-1}$. Here~\eqref{rcond} becomes
\be{rc2}
(1+\del^2)G(\del+\del^3)\le\del,
\ee{rc2}
and observing that $g'(0)=0$, we see that $G(u)=O(u^2)$ as $u\to 0$. Thus the inequality~\eqref{rc2} holds
for all sufficiently small $\del$, that is,
for all sufficiently large $\lam$, as desired.

Finally, the bound~\eqref{sigb} on the coefficients follows from the integral formula
$$
\sigma_n=\frac{\sigma^{(n)}(0)}{n!}=\frac{1}{2\pi i}\int_{|s|=1}\frac{\sigma(s)-s}{s^{n+1}}\:ds
$$
for $n\ge 2$, using~\eqref{qbn}.~\fp
\vv

\noindent {\bf Proof of Theorem~5.2.}
We have a unique (up to positive scalar multiple) positive periodic solution $x(t)$
with eigenvalue $\kappa$, by Theorem~5.1. Also, we note that
$$
\int^t_{t-r(t)}\rho(s)\:ds=r(t)
$$
and thus by~\eqref{radbnd} we have the bounds
\be{kap}
0<2\pi m-\lam+1\le\kappa\le2\pi m+\lam-1
\ee{kap}
for the eigenvalue.
The differential equation~\eqref{u++} takes the form
$$
\kappa\dot x(t)=x(t)-\dot\eta(t)x(\eta(t)),
$$
with $\eta(t)$ as in~\eqref{eta}. We now apply the conjugacy
of Lemma~4.1 at the point $t_0=0$. Letting $t=\sigma(\wt t)$ and $y(\wt t)=x(\sigma(\wt t))$ for $\wt t$ near zero,
we obtain the equation
$$
\kappa\dot y(\wt t)=\dot\sigma(\wt t)y(\wt t)-\dot\sigma(\wt t)\dot\eta(\sigma(\wt t))y(\lambda\wt t),
$$
and noting that $\eta(\sig(\wt t))=\sig(\lam \wt t)$ and hence
$\dot\sigma(\wt t)\dot\eta(\sigma(\wt t))=\lambda\dot\sigma(\lambda\wt t)$, we rewrite this equation as
$$
\kappa\dot y(\wt t)=\dot\sigma(\wt t)y(\wt t)-\lambda\dot\sigma(\lambda\wt t)y(\lambda\wt t).
$$
With the Taylor series of $\sigma(\wt t)$ given by~\eqref{sg}, we have
$$
\alpha_n=\frac{(n+1)\sigma_{n+1}}{\kappa},\qquad
\beta_n=-\frac{(n+1)\lambda^{n+1}\sigma_{n+1}}{\kappa},
$$
for $n\ge 0$, for the Taylor coefficients of $\alpha(\wt t)$ and $\beta(\wt t)$ as in~\eqref{yeqq}, \eqref{abg}.
Using the above formulas for $\alpha_n$ and $\beta_n$, we proceed as in equation~\eqref{yr} and write
a recursion formula for $y_{n+1}$ in terms of $y_k$ for $0\le k\le n$ and of $\sigma_k$ for $1\le k\le n+1$.
Equation~\eqref{cy} thus takes the form $w_k=\xi_ky_k$ where
$$
\xi_k=\frac{(-1)^k\kappa^kk!}{\lambda^{k(k+1)/2}}.
$$
(The quantity $\xi_k$ is related to, but slightly different from, the quantity
$\theta_k$ in~\eqref{the}.)
After substituting for $y_k$ in terms of $w_k$ in the formula for $y_{n+1}$, one sees after some manipulation
that the resulting formula~\eqref{cr} for $w_{n+1}$ in the proof of Theorem~4.2 takes the form
\be{crc}
w_{n+1}=\bigg(1-\frac{1}{\lambda^{n+1}}\bigg)\bigg(w_n
+\sum_{k=0}^{n-1}\frac{\xi_n}{\xi_k}
\bigg((n-k+1)\sigma_{n-k+1}\bigg)w_k\bigg).
\ee{crc}
We now obtain estimates for the quantities $w_n$ and thereby for $w_\infty$ much as
in the proof of Theorem~4.2, but with some differences. Let us first observe that
$$
\kappa\le 2\pi m+\lam-1\le 3\lam,
$$
by~\eqref{kap} and the choice of $m$ as in the statement of the theorem. It follows that
$$
\bigg|\frac{\xi_k}{\xi_{k-1}}\bigg|=\frac{\kappa k}{\lambda^k}\le\frac{3k}{\lambda^{k-1}}\le
\left\{\begin{array}{ll}
3, &\bb \quad\hbox{for }k=1,\\
\\
\ds{\frac{1}{2},} &\bb \quad\hbox{for }k>1,
\end{array}\right.
$$
where we assume that $\lambda$ is large enough that $6k\le\lambda^{k-1}$
for every $k>1$. (More precisely, $\lam\ge 12$ suffices here.) Therefore,
$$
\bigg|\frac{\xi_n}{\xi_k}\bigg|
=\bigg|\frac{\xi_n}{\xi_{n-1}}\bigg|\bigg|\frac{\xi_{n-1}}{\xi_k}\bigg|
\le\bigg(       \frac{3n}{\lambda^{n-1}}\bigg)   \bigg(\frac{3}{2^{n-k-2}}\bigg)     \le 18\bigg(\frac{2^k}{\lam^{n-1}}\bigg)
$$
for $0\le k\le n-1$. It now follows from this, from~\eqref{crc}, and from the bound~\eqref{sigb} in Lemma~5.3, that
$$
\begin{array}{lcl}
|w_{n+1}-w_n| &\bb \le &\bb
\ds{\frac{|w_n|}{\lambda^{n+1}}
+18\sum_{k=0}^{n-1}\bigg(\frac{2^k}{\lam^n}\bigg)(n-k+1)|w_k|}\\
\\
&\bb \le &\bb
\ds{\bigg(\frac{2}{\lam}\bigg)^n\bigg(1+18\sum_{j=1}^\infty\frac{j+1}{2^j}\bigg)\max_{0\le k\le n}|w_k|
=\bigg(\frac{2}{\lam}\bigg)^nK\max_{0\le k\le n}|w_k|,}
\end{array}
$$
for every $n\ge 0$,
with the above equality defining the constant $K$. In fact, for $n=0$ a stronger bound is available, as
one sees from~\eqref{crc} that
$$
|w_1-w_0|=\frac{|w_0|}{\lam},
$$
and so in any case we have
\be{cdf}
|w_{n+1}-w_n|\le H_n(\lam)\max_{0\le k\le n}|w_k|,\qquad\mbox{where}\qquad
H_n(\lam)=\left\{\begin{array}{ll}
\ds{\frac{1}{\lam},} &\bb \quad\mbox{for }n=0,\\
\\
\ds{\bigg(\frac{2}{\lam}\bigg)^nK,} &\bb \quad\mbox{for }n>0.
\end{array}\right.
\ee{cdf}
It follows, much as in the proof of Theorem~4.2, that
$$
|w_n|\le\bigg(\prod_{k=0}^{n-1}(1+H_k(\lambda))\bigg)|w_0|\le\bigg(\prod_{k=0}^\infty(1+H_k(\lambda))\bigg)|w_0|,
$$
and so using this bound in~\eqref{cdf} gives
$$
|w_{n+1}-w_n|\le H_n(\lambda)\bigg(\prod_{k=0}^\infty(1+H_k(\lambda))\bigg)|w_0|.
$$
Upon summing we obtain
\be{om}
|w_\infty-w_0|\le\Omega(\lambda)|w_0|,\qquad
\Omega(\lambda)=\bigg(\sum_{n=0}^\infty H_n(\lambda)\bigg)\bigg(\prod_{k=0}^\infty(1+H_k(\lambda))\bigg),
\ee{om}
with the above formula serving as the definition of $\Omega(\lambda)$. We note in particular that the infinite
sum and product in this formula have finite positive values as long as $\lambda>2$.

If $\Omega(\lambda)<1$ then our desired conclusion, that
$w_\infty\ne 0$ if $w_0\ne 0$, follows from~\eqref{om}. Indeed, this
is the case for all large $\lambda$, as follows directly from the fact that
$\ds{\lim_{\lambda\to\infty}H_n(\lambda)=0}$ monotonically for every $n\ge 0$.
With this, the theorem is proved.~\fp
\vv

We end this section with the following result, which shows that
for certain $\lam$ and $m$, the periodic solution considered in Theorem~5.2
is analytic for certain values of $t$. More precisely, in addition to the
expansive fixed point $t_0=0$ at which, by Theorem~5.2, analyticity fails, there is also
a contractive fixed point $t_{00}$ in a neighborhood of which the solution is 
analytic. Thus one has both $\mathcal{N}\ne\emptyset$ and $\mathcal{A}\ne\emptyset$.
\vv

\noindent {\bf Theorem~5.4.} \bxx
Consider the integral equation~\eqref{ie} where $3<\lam<2\pi m+1$, with $m$ an integer.
Suppose further that there exists an integer $n$ satisfying
\be{pq}
\bigg((\lam-1)^2-4\bigg)^{1/2}<2\pi n<\lam-1.
\ee{pq}
Then there exists some $t_{00}\in\IR$ such that the unique positive $2\pi$-periodic
solution $x(t)$ of~\eqref{ie} is analytic for $t$ in some neighborhood of $t_{00}$.
More precisely, one has that $\eta(t_{00})=t_{00}$ and $|\dot\eta(t_{00})|<1$ where
$\eta(t)=t-r(t)+2\pi(m-n)$, as in Theorem~2.1.
\exx
\vv

\noindent {\bf Proof.}
Let
$$
t_{00}=\frac{\pi}{2}+\tau,\qquad
\tau=\arccos\bigg(\frac{2\pi n}{\lam-1}\bigg),
$$
with $0<\tau<\frac{\pi}{2}$.
Then noting from~\eqref{ie} and from the statement of the present proposition that
$\eta(t)=t+(\lam-1)\sin t-2\pi n$, it follows that
$$
\eta(t_{00})=t_{00}+(\lam-1)\sin t_{00}-2\pi n=t_{00}+(\lam-1)\cos\tau-2\pi n=t_{00}.
$$
Further,
$$
\dot\eta(t_{00})=1+(\lam-1)\cos t_{00}=1-(\lam-1)\sin\tau=1-\bigg((\lam-1)^2-(2\pi n)^2\bigg)^{1/2},
$$
and one sees immediately from~\eqref{pq} that $|\dot\eta(t_{00})|<1$. The desired conclusions
now follow directly from Theorem~2.1.~\fp

\section{Open Questions}

Many open questions remain. For a given solution of interest,
very generally one wishes to determine the set $\mathcal{A}$ of analyticity and its complement $\mathcal{N}$.
Toward this end, one open problem is to extend the results on expansive fixed points, such as Theorems~4.2 and~4.5,
to the case of systems (that is, with $x\in\IR^N$) and to the case of expansive periodic points (that is,
with $\eta^M(t_0)=t_0$ and $|\dot\eta^M(t_0)|>1$); and further, to extend such results
from linear to general nonlinear equations~\eqref{nl}.
As noted in one of the remarks following the proof of Theorem~4.2,
the proof of that result does not carry over to linear systems due to commutativity issues
with the matrix-valued coefficients.
In dealing with periodic points of $\eta$, one naturally encounters
the system~\eqref{nl2} as in the proof of Theorem~2.1, so the commutativity problem arises there also.

Another issue is that even in the case of an expansive fixed point in a scalar linear system, if $w_\infty=0$ then
by Theorem~4.4 there is no assurance that a solution of interest is analytic at the point in question.
One might ask what other conditions are sufficient to give analyticity under such circumstances.

Questions remain in the case of points $t_0$ which are neither periodic for the map $\eta$,
nor are in the basin of attraction of a periodic point. This includes
time-periodic systems for which $\eta$ is a homeomorphism of the circle $S^1$ with an irrational
rotation number, but it also includes more general cases in which the orbit $\eta^n(t_0)$
of some $t_0$ exhibits a chaotic character. As noted at the end of Section~3,
it would be interesting, even for linear equations with periodic coefficients, and where
$\eta:S^1\to S^1$ is a homeomorphism with an irrational
rotation number, to find an example with a nowhere analytic but everywhere $C^\infty$
periodic solution. A characterization of the set of rotation numbers for which this
is possible is also of interest. One could also consider the case of quasiperiodic and almost periodic systems,
and it would not be unreasonable to expect issues involving small divisors to occur in their analysis.
Additionally, one might consider analyticity properties of solutions on the unstable manifold of an
equilibrium or periodic orbit.

A very basic problem is to extend the results herein to equations with multiple delays, such
as equation~\eqref{mult}, where $r_k=r_k(t)$ for $1\le k\le M$ are given analytic functions of $t$.
(Of course the nonlinearity $f$ is also analytic.) A related issue is the paucity of results
on existence of interesting solutions (for example, periodic solutions) in equations with
multiple delays, although there typically exist Floquet solutions for linear equations of the form~\eqref{mult}.

One might also consider more general classes of equations, for example integral equations such as
\be{inteq2}
x(t)=\int^t_{t-r(t)}f(t,s,x(s))\:ds.
\ee{inteq2}
In contrast to the integral equation~\eqref{inteq} considered in Section~5, differentiating~\eqref{inteq2}
does not lead to a differential equation to which our results can readily be applied.
Another interesting problem is to obtain analyticity results for equations with almost constant delays
$r(t)=r_0+\eps r_1(t)$. A perturbation analysis in this spirit would be
quite natural for equation~\eqref{inteq} with the periodicity conditions~\eqref{v1},
with $r_0=2\pi m$ and $r_1$ of period $2\pi$, where we note that for the unperturbed problem with $\eps=0$,
the eigenfunction would be a constant $x(t)\equiv 1$.

A still unsolved question raised in~\cite{nkr} concerns the equation
$$
\kappa x(t)=\int_{t-r}^t\rho(s)x(s)\:ds
$$
where here $r>0$ is a given constant and $\rho:\IR\to(0,\infty)$ is analytic with $\rho(t+2\pi)=\rho(t)$
for all $t\in\IR$. Let $\kappa=\kappa(r)$ and $x(t)=x(t,r)=x(t+2\pi,r)>0$ denote the eigenvalue and eigenfunction
given by Theorem~5.1, normalized so that, say, $x(0,r)=1$. It has been proved in~\cite{nkr} that $\kappa(r)$
is a $C^\infty$ function of $r$, and the results in~\cite{nuss} for problems with constant delays imply that
$x(t,r)$ is analytic in $t$ for every fixed $r$. However, it is not known whether $\kappa(r)$ is analytic
in $r$ for any range of $r$, and similarly for $x(t,r)$.

A very significant problem, and which in a sense is partly the motivation for the present study,
is to understand analyticity properties of solutions of equations with state-dependent delays. One of the
simplest such systems would take the form~\eqref{state}, although many other methods of incorporating
state-dependent delays (for example, implicitly defined delays) are possible. Note that the ``linear'' version
of this equation, namely equation~\eqref{slin} with $\eta(t)$ replaced with $t-r(x(t))$, is of course
not linear at all in general. In any case, if $x(t)$ is a solution of such a state-dependent problem
with $x(t_0)=x_0$ and with $t_0-r(x_0)=t_0$ holding, that is, $r(x_0)=0$, then one could still distinguish
contractive from expansive fixed points via the magnitude of $|\frac{d}{dt}(t-r(x(t)))|=|1-r'(x(t))\dot x(t)|$
at $t=t_0$, and one could still expand $x(t)$ about this point to obtain uniquely determined Taylor coefficients $x_n$.
However, the implementation of the Hartman--Grobman result Lemma~4.1 is not clear, as one does not know {\it a priori}
whether the map $t\to t-r(x(t))$ is analytic; and the resulting recursion relation
for the coefficients $x_n$, while well-defined, would
be much more complicated than~\eqref{yr}.

A somewhat technical problem is to extend Theorem~4.2 to the case when $\bet(0)=\bet_0=0$. Even with simple
examples, such as
$$
\dot y(t)=t^my(\lam t),\qquad y(0)=1,
$$
where $m>0$ and $|\lam|>1$, one encounters new issues.
In particular, from the recursion~\eqref{yr} it follows that $y_n=0$ if $n\not\equiv 0$ (mod $m+1$),
while for $n=(m+1)k$ one has
$$
y_{(m+1)k}=\frac{\lam^{(m+1)k(k-1)/2}}{(m+1)^kk!}.
$$
One sees that $\ds{\lim_{n\to\infty}}|y_n|^{1/n}$ does not exist, although the lim~sup of this quantity is $\infty$,
and so no analytic solution exists.

Another technical issue is to obtain sharp conditions on $f$, in Theorem~3.4, in the system case $N>1$,
as noted in the remarks following the proof of that result. One requires
some sort of nontrivial dependence of $f$ on the delay term $v$; however, such conditions, as condition~(1) of that theorem
for the scalar case, should be broad enough to allow for a nonlinear dependence
on the delay term more general than the linear dependence of condition~(2).
Quite possibly, algebraic techniques from ring theory in the spirit of Neelon's work~\cite{neelon} would be helpful.

Many questions involving the quantity $w_\infty$ in Theorem~4.2 present themselves. One can regard $w_\infty$ as a function
of the coefficients $\alp$, $\bet$, $\gam$, of $\eta$, and of the initial condition $y_0$, and we may ask how $w_\infty$
varies with respect to these. Some pertinent questions are: Does $w_\infty$ vary continuously with respect to these data, in an appropriate norm?
Does it vary smoothly, or even analytically? Is it the case that generically
(namely for a residual set of data $\alp$, $\bet$, $\gam$, $\eta$, and $y_0$) one has $w_\infty\ne 0$, or more simply,
does $w_\infty$ have a nontrivial dependence on these data? (See in particular the question posed in the remark following
the proof of Theorem~4.2.)
Studying specific simple systems such as
$$
\dot y(t)=(1+zt)y(\lam t),\qquad y(0)=1,
$$
where $z\in\IR$ and $\lam\in(-\infty,1)\cup(1,\infty)$,
might provide some insights.
Here one has $w_\infty=w_\infty(z,\lam)$ as a function of $z$ and $\lam$,
and both theoretical and numerical studies could be undertaken.

More broadly, many questions on analyticity remain even for problems with constant delays.
Suppose that $x(t)$ is a periodic solution of a system~\eqref{1a}
where $f:\IR^{N(M+1)}\to\IR^N$ is analytic and the delays $r_k\ge 0$ for $1\le k\le M$ are given constants;
one may think of well-known systems such as Wright's equation
$$
\dot x(t)=-\alp x(t-1)(1+x(t)),
$$
or the Mackey--Glass equation (with analytic $g$)
$$
\alp\dot x(t)=-x(t)+g(x(t-1)),
$$
which have periodic solutions for a range of the parameter $\alp$.
Then as noted, from~\cite{nuss} the solution $x(t)$ is analytic, at least for real $t$. However,
virtually nothing is known about the analytic continuation of this solution into the complex plane,
that is, for complex $t$.

Analyticity issues arise in parameterized systems with constant delays. Consider such an autonomous system
$$
\dot x(t)=f(x(t),x(t-r_1),x(t-r_2),\ldots,x(t-r_M),\alp),
$$
with the nonlinearity $f:\IR^{N(M+1)}\times\IR\to\IR^N$ analytic
in all its arguments, and suppose for some parameter value $\alp=\alp_0$ there is
a periodic solution $x_0(t)$ of period $p_0$. Suppose further
that this periodic solution is hyperbolic, that is, all of its nontrivial characteristic
multipliers $\mu$ satisfy $|\mu|\ne 1$. Then for all $\alp$ near $\alp_0$ one has a nearby
periodic solution $x(t,\alp)$ of period $p(\alp)$ near $p_0$, obtained in a
standard fashion as the Poincar\'e continuation of $x_0(t)$. Certainly $x(t,\alp)$
is analytic in $t$ for each fixed $\alp$.
However, whether or not $x(t,\alp)$ and $p(\alp)$
depend analytically on the parameter $\alp$ is not known.

\end{document}